\documentclass[a4paper]{article}

\usepackage{MyLaTeX,MyArticle}
\usepackage{lscape}
\usepackage{hyperref}
\RequirePackage{pgf,pgffor}
\usepackage{tikz}
\usetikzlibrary{matrix,arrows,decorations.pathmorphing}

\numberwithin{equation}{section}

\title{An Atlas of Modular Representation Theory\\Version 1: Information on $\Ext^1$ for simple modules for groups of Lie type in defining characteristic over small fields}
\author{David A. Craven}
\date{\today}
\begin{document}
\maketitle

\section{Introduction}

This document is the first iteration of an attempt to collate information about small-rank groups of Lie type over small fields, and their representation theory over the defining field. This information is important in the author's work on subgroup structure of exceptional groups of Lie type.

The most important information in that work is information about $\Ext^1$ between simple modules, and so in Version 1 of this document, that data is almost all of the data available. In addition, a lot of information about the dimensions of the simple and Weyl modules is included.

More generally, one may expect to include details about the socle structure of the projective modules, Jordan block structure of the action of unipotent elements, decompositions of symmetric and exterior powers of simple modules, and tensor products of modules, traces of semisimple elements and so on.

The ideal place for such information is a dedicated website, connected to a database that could be queried to produce the information required. This document, while imperfect, will have to suffice for now.

\medskip

Note that (almost) all of these data were constructed via computer and with Magma. There are various articles that compute $\Ext^1$ between simple modules for low-rank Lie type groups over small fields, but given that long computational papers are notoriously error prone, the author felt it best to compute the table here independently, so that we now have two sources for (some of) these tables.

The labelling conventions for our simple modules follow the author's recent papers on maximal subgroups, with a dictionary between this and highest-weight labellings being given at the start of each section.

\newpage

\tableofcontents

\newpage
\section{\texorpdfstring{$\SL_3(2)$}{SL(3,2)}}

\subsection*{Block information}

\begin{center}\begin{tabular}{ccc}
\hline Block & Kernel & Modules
\\\hline Block 1 & 1 & $1$, $3$, $3^*$
\\Block 2 & 1 & $8$
\\\hline\end{tabular}\end{center}

\subsection*{Highest weight labels}

\begin{center}\begin{tabular}{ccccc}
\hline $\lambda$&00&10&01&11
\\\hline $L(\lambda)$ &$1$&$3$&$3^*$&$8$\\ $\dim(W(\lambda))$&$1$&$3$&$3$&$8$
\\\hline
\end{tabular}\end{center}

\subsection*{Extension information}

\begin{center}\begin{tabular}{c|ccc}
Block 1&$1$&$3$&$3^*$
\\\hline $1$&$0$&$1$&$1$
\\ $3$&$1$&$0$&$1$
\\ $3^*$&$1$&$1$&$0$
\\\end{tabular}\end{center}

\newpage
\section{\texorpdfstring{$\SL_3(3)$}{SL(3,3)}}

\subsection*{Block information}

\begin{center}\begin{tabular}{ccc}
\hline Block & Kernel & Modules
\\\hline Block 1 & 1 & $1$, $3$, $3^*$, $6$, $6^*$, $7$, $15$, $15^*$
\\Block 2 & 1 & $27$
\\\hline\end{tabular}\end{center}

\subsection*{Identification of modules}

\begin{enumerate}
\item Fix $3$ and $3^*$.
\item $6=S^2(3^*)$.
\item $15=\Lambda^2(6^*)$.
\end{enumerate}

\subsection*{Highest weight labels}

\begin{center}\begin{tabular}{cccccccccc}
\hline $\lambda$&00&10&01&02&20&11&21&12&22
\\\hline $L(\lambda)$ &$1$&$3$&$3^*$&$6$&$6^*$&$7$&$15$&$15^*$&$27$\\ $\dim(W(\lambda))$&$1$&$3$&$3$&$6$&$6$&$8$&$15$&$15$&$27$
\\\hline
\end{tabular}\end{center}

\subsection*{Extension information}

\begin{center}\begin{tabular}{c|cccccccc}
Block 1&$1$&$3$&$3^*$&$6$&$6^*$&$7$&$15$&$15^*$
\\\hline $1$&$0$&$0$&$0$&$0$&$0$&$1$&$1$&$1$
\\ $3$&$0$&$0$&$0$&$0$&$1$&$1$&$0$&$0$
\\ $3^*$&$0$&$0$&$0$&$1$&$0$&$1$&$0$&$0$
\\ $6$&$0$&$0$&$1$&$0$&$1$&$0$&$0$&$0$
\\ $6^*$&$0$&$1$&$0$&$1$&$0$&$0$&$0$&$0$
\\ $7$&$1$&$1$&$1$&$0$&$0$&$0$&$0$&$0$
\\ $15$&$1$&$0$&$0$&$0$&$0$&$0$&$0$&$0$
\\ $15^*$&$1$&$0$&$0$&$0$&$0$&$0$&$0$&$0$
\end{tabular}\end{center}

\newpage
\section{\texorpdfstring{$\SL_3(4)$}{SL(3,4)}}

\subsection*{Block information}

\begin{center}\begin{tabular}{ccc}
\hline Block & Kernel & Modules
\\ \hline 1 & 3 & $8_1$, $8_2$, $9_{1,2}$, $9_{1,2}^*$
\\ 2 & 3 & $64$
\\ 3 & 1 & $3_1$, $3_2^*$, $\bar 9_{1,2}^*$ $24_{1,2}$, $24_{2,1}^*$
\\ 4 & 1 & Dual to block 3
\\\hline\end{tabular}\end{center}

\subsection*{Identification of modules}

We write $3_2$ for the twist of $3_1$, so that that $3_1$ and $3_2$ restrict to the same module for $\SL_3(2)$. Write $8_1$ for the module such that $8_1\otimes 3_2$ is irreducible.
\begin{enumerate}
\item $9_{1,2}=3_1\otimes 3_2$
\item $\bar 9_{1,2}=3_1\otimes 3_2^*$
\item $24_{1,2}=3_1\otimes 8_2$ and $24_{2,1}=3_2\otimes 8_1$
\end{enumerate}

\subsection*{Extension information}

\begin{center}\begin{tabular}{c|ccccc}
Block 1 &$1$&$8_1$&$8_2$&$9_{1,2}$&$9_{1,2}^*$
\\\hline $1$&$0$&$0$&$0$&$2$&$2$
\\ $8_1$&$0$&$0$&$0$&$1$&$1$
\\ $8_2$&$0$&$0$&$0$&$1$&$1$
\\ $9_{1,2}$&$2$&$1$&$1$&$0$&$0$
\\ $9_{1,2}^*$&$2$&$1$&$1$&$0$&$0$
\end{tabular}\end{center}

\begin{center}\begin{tabular}{c|ccccc}
Block 3 &$3_1$&$3_2^*$&$\bar 9_{1,2}^*$&$24_{1,2}$&$24_{2,1}^*$
\\\hline $3_1$&$0$&$2$&$1$&$0$&$1$
\\ $3_2^*$&$2$&$0$&$1$&$1$&$0$
\\ $\bar 9_{1,2}^*$&$1$&$1$&$0$&$0$&$0$
\\ $24_{1,2}$&$0$&$1$&$0$&$0$&$0$
\\ $24_{2,1}^*$&$1$&$0$&$0$&$0$&$0$
\end{tabular}\end{center}

\newpage
\section{\texorpdfstring{$\SL_3(5)$}{SL(3,5)}}

\subsection*{Block information}

\begin{center}\begin{tabular}{ccc}
\hline Block & Kernel & Modules
\\\hline Block 1 & $1$ & $1$, $3$, $3^*$, $6$, $6^*$, $8$, $10$, $10^*$, $15_1$, $15_1^*$, $15_2$, $15_2^*$
\\ && $18$, $18^*$, $19$, $35$, $35^*$, $39$, $39^*$, $60$, $60^*$, $63$, $90$, $90^*$
\\Block 2 & $1$ & $125$
\\\hline\end{tabular}\end{center}

\subsection*{Identification of modules}

\begin{enumerate}
\item Fix $3$ and $3^*$.
\item $6=S^2(3^*)$.
\item $10$ is a submodule of $3^*\otimes 6$.
\item $15_1$ is a submodule of $S^2(6^*)$.
\item $15_2=\Lambda^2(6^*)$.
\item $18$ is a submodule of $6\otimes 10$.
\item $35$ is a submodule of $\Lambda^2(10^*)$.
\item $39$ is a submodule of $\Lambda^2(15_1^*)$.
\item $60$ is a submodule of $6\otimes 19$.
\item $90$ is a submodule of $\Lambda^2(18^*)$.
\end{enumerate}

\subsection*{Highest weight labels}

\begin{center}\begin{tabular}{ccccccccccccccc}
\hline $\lambda$&00&10&01&02&20&11&03&30&40&04&21&12&13&31
\\\hline $L(\lambda)$ &$1$&$3$&$3^*$&$6$&$6^*$&$8$&$10$&$10^*$&$15_1$&$15_1^*$&$15_2$&$15_2^*$&$18$&$18^*$
\\ $\dim(W(\lambda))$&$1$&$3$&$3$&$6$&$6$&$8$&$10$&$10$&$15$&$15$&$15$&$15$&$24$&$24$
\\\hline
\end{tabular}\end{center}

\begin{center}\begin{tabular}{cccccccccccc}
\hline $\lambda$&22&41&14&32&23&24&42&33&43&34&44
\\\hline $L(\lambda)$ &$19$&$35$&$35^*$&$39$&$39^*$&$60$&$60^*$&$63$&$90$&$90^*$&$125$
\\ $\dim(W(\lambda))$&$27$&$35$&$35$&$42$&$42$&$60$&$60$&$64$&$90$&$90$&$125$
\\\hline
\end{tabular}\end{center}

\subsection*{Extension information}

\begin{center}\begin{footnotesize}\begin{tabular}{c|ccccccccccccccccccccc}
&$1$&$3$&$3^*$&$6$&$6^*$&$8$&$10$&$10^*$&$15_1$&$15_1^*$&
$15_2$&$15_2^*$&$18$&$18^*$&$19$&$35$&$35^*$&$39$&$60$&$63$&$90$
\\\hline $1$&$0$&$0$&$0$&$0$&$0$&$0$&$0$&$0$&$0$&$0$&$0$&$0$&$0$&$0$&$0$&$0$&$0$&$1$&$0$&$1$&$0$
\\ $3$&$0$&$0$&$0$&$0$&$0$&$0$&$0$&$0$&$0$&$0$&$0$&$0$&$0$&$1$&$1$&$0$&$0$&$1$&$0$&$1$&$0$
\\ $3^*$&$0$&$0$&$0$&$0$&$0$&$0$&$0$&$0$&$0$&$0$&$0$&$0$&$1$&$0$&$1$&$0$&$0$&$1$&$0$&$1$&$0$
\\ $6$&$0$&$0$&$0$&$0$&$0$&$0$&$0$&$0$&$0$&$0$&$0$&$0$&$1$&$0$&$1$&$1$&$0$&$0$&$0$&$0$&$0$
\\ $6^*$&$0$&$0$&$0$&$0$&$0$&$0$&$0$&$0$&$0$&$0$&$0$&$0$&$0$&$1$&$1$&$0$&$1$&$1$&$0$&$0$&$1$
\\ $8$&$0$&$0$&$0$&$0$&$0$&$0$&$0$&$0$&$0$&$0$&$0$&$0$&$1$&$1$&$1$&$0$&$0$&$1$&$1$&$0$&$0$
\\ $10$&$0$&$0$&$0$&$0$&$0$&$0$&$0$&$0$&$0$&$0$&$0$&$0$&$0$&$1$&$0$&$0$&$0$&$0$&$1$&$0$&$0$
\\ $10^*$&$0$&$0$&$0$&$0$&$0$&$0$&$0$&$0$&$0$&$0$&$0$&$0$&$1$&$0$&$0$&$0$&$0$&$0$&$0$&$0$&$0$
\\ $15_1$&$0$&$0$&$0$&$0$&$0$&$0$&$0$&$0$&$0$&$1$&$0$&$1$&$0$&$0$&$0$&$0$&$0$&$0$&$0$&$0$&$0$
\\ $15_1^*$&$0$&$0$&$0$&$0$&$0$&$0$&$0$&$0$&$1$&$0$&$1$&$0$&$0$&$0$&$0$&$0$&$0$&$0$&$0$&$0$&$0$
\\ $15_2$&$0$&$0$&$0$&$0$&$0$&$0$&$0$&$0$&$0$&$1$&$0$&$0$&$0$&$0$&$0$&$1$&$0$&$0$&$0$&$1$&$0$
\\ $15_2^*$&$0$&$0$&$0$&$0$&$0$&$0$&$0$&$0$&$1$&$0$&$0$&$0$&$0$&$0$&$0$&$0$&$1$&$1$&$0$&$1$&$0$
\\ $18$&$0$&$0$&$1$&$1$&$0$&$1$&$0$&$1$&$0$&$0$&$0$&$0$&$0$&$0$&$0$&$0$&$0$&$0$&$0$&$0$&$0$
\\ $18^*$&$0$&$1$&$0$&$0$&$1$&$1$&$1$&$0$&$0$&$0$&$0$&$0$&$0$&$0$&$0$&$0$&$0$&$0$&$0$&$0$&$0$
\\ $19$&$0$&$1$&$1$&$1$&$1$&$1$&$0$&$0$&$0$&$0$&$0$&$0$&$0$&$0$&$0$&$0$&$0$&$0$&$0$&$0$&$0$
\\ $35$&$0$&$0$&$0$&$1$&$0$&$0$&$0$&$0$&$0$&$0$&$1$&$0$&$0$&$0$&$0$&$0$&$0$&$0$&$0$&$0$&$0$
\\ $35^*$&$0$&$0$&$0$&$0$&$1$&$0$&$0$&$0$&$0$&$0$&$0$&$1$&$0$&$0$&$0$&$0$&$0$&$0$&$0$&$0$&$0$
\\ $39$&$1$&$1$&$1$&$0$&$1$&$1$&$0$&$0$&$0$&$0$&$0$&$1$&$0$&$0$&$0$&$0$&$0$&$0$&$0$&$0$&$0$
\\ $39^*$&$1$&$1$&$1$&$1$&$0$&$1$&$0$&$0$&$0$&$0$&$1$&$0$&$0$&$0$&$0$&$0$&$0$&$0$&$0$&$0$&$0$
\\ $60$&$0$&$0$&$0$&$0$&$0$&$1$&$1$&$0$&$0$&$0$&$0$&$0$&$0$&$0$&$0$&$0$&$0$&$0$&$0$&$0$&$0$
\\ $60^*$&$0$&$0$&$0$&$0$&$0$&$1$&$0$&$1$&$0$&$0$&$0$&$0$&$0$&$0$&$0$&$0$&$0$&$0$&$0$&$0$&$0$
\\ $63$&$1$&$1$&$1$&$0$&$0$&$0$&$0$&$0$&$0$&$0$&$1$&$1$&$0$&$0$&$0$&$0$&$0$&$0$&$0$&$0$&$0$
\\ $90$&$0$&$0$&$0$&$0$&$1$&$0$&$0$&$0$&$0$&$0$&$0$&$0$&$0$&$0$&$0$&$0$&$0$&$0$&$0$&$0$&$0$
\\ $90^*$&$0$&$0$&$0$&$1$&$0$&$0$&$0$&$0$&$0$&$0$&$0$&$0$&$0$&$0$&$0$&$0$&$0$&$0$&$0$&$0$&$0$
\end{tabular}\end{footnotesize}\end{center}

\newpage
\section{\texorpdfstring{$\SL_3(7)$}{SL(3,7)}}

\subsection*{Block information}

\begin{center}\begin{tabular}{ccc}
\hline Block & Kernel & Modules
\\\hline Block 1 & $3$ & $1$, $8$, $10$, $10^*$, $27$, $28$, $28^*$, $35$, $35^*$, $37$, $71$, $71^*$, $117$, $154$, $154^*$, $215$
\\Block 2 & $3$ & $343$
\\ Block 3 & $1$ & $3$, $6$, $15_1$, $15_2$, $21$, $24$, $33$, $36$, $42$, $63$, $75$, $105$, $114$, $162$, $210$, $273$
\\ Block 4 & $1$ & Dual to block 3
\\\hline\end{tabular}\end{center}

\subsection*{Identification of modules}

\begin{enumerate}
\item Fix $3$ and $3^*$. This fixes all modules in Blocks 3 and 4 but $15_1$ and $15_2$.
\item $10$ is a submodule of $3^*\otimes 6$.
\item $15_1$ is a submodule of $S^2(6^*)$.
\item $15_2=\Lambda^2(6^*)$.
\item $28$ is a submodule of $S^2(10^*)$.
\item $35$ is a submodule of $\Lambda^2(10^*)$.
\item $71$ is a submodule of $8\otimes 28^*$.
\item $154$ is a submodule of $10\otimes 37$.
\end{enumerate}

\subsection*{Highest weight labels}

\begin{center}\begin{tabular}{cccccccccccccccccc}
\hline $\lambda$&00&11&03&30&22&60&06&41&14&33&25&52&44&36&63&55&66
\\\hline $L(\lambda)$ &$1$&$8$&$10$&$10^*$&$27$&$28$&$28^*$&$35$&$35^*$&$37$&$71$&$71^*$&$117$&$154$&$154^*$&$215$&$343$
\\ $\dim(W(\lambda))$&$1$&$8$&$10$&$10$&$27$&$28$&$28$&$35$&$35$&$64$&$81$&$81$&$125$&$154$&$154$&$216$&$343$
\\\hline
\end{tabular}\end{center}

\begin{center}\begin{tabular}{ccccccccccccccccc}
\hline $\lambda$&10&02&40&21&05&13&51&24&32&16&43&62&35&54&46&65
\\\hline $L(\lambda)$ &$3$&$6$&$15_1$&$15_2$&$21$&$24$&$33$&$36$&$42$&$63$&$75$&$105$&$114$&$162$&$210$&$273$
\\ $\dim(W(\lambda))$&$3$&$6$&$15$&$15$&$21$&$24$&$48$&$60$&$42$&$63$&$90$&$105$&$120$&$165$&$210$&$273$
\\\hline
\end{tabular}\end{center}

\subsection*{Extension information}

\begin{center}\begin{tabular}{c|cccccccccccccccc}
Block 1&$1$&$8$&$10$&$10^*$&$27$&$28$&$28^*$&$35$&$35^*$&$37$&$71$&$71^*$&$117$&$154$&$154^*$&$215$
\\ \hline $1$&$0$&$0$&$0$&$0$&$0$&$0$&$0$&$0$&$0$&$0$&$1$&$1$&$0$&$0$&$0$&$1$
\\ $8$&$0$&$0$&$0$&$0$&$0$&$0$&$0$&$0$&$0$&$2$&$1$&$1$&$3$&$0$&$0$&$0$
\\ $10$&$0$&$0$&$0$&$0$&$0$&$0$&$0$&$0$&$0$&$1$&$1$&$1$&$1$&$0$&$0$&$1$
\\ $10^*$&$0$&$0$&$0$&$0$&$0$&$0$&$0$&$0$&$0$&$1$&$1$&$1$&$1$&$0$&$0$&$1$
\\$27$&$0$&$0$&$0$&$0$&$0$&$0$&$0$&$0$&$0$&$1$&$1$&$1$&$2$&$1$&$1$&$0$
\\$28$&$0$&$0$&$0$&$0$&$0$&$0$&$1$&$0$&$1$&$0$&$0$&$0$&$0$&$0$&$0$&$0$
\\$28^*$&$0$&$0$&$0$&$0$&$0$&$1$&$0$&$1$&$0$&$0$&$0$&$0$&$0$&$0$&$0$&$0$
\\$35$&$0$&$0$&$0$&$0$&$0$&$0$&$1$&$0$&$0$&$0$&$0$&$1$&$0$&$0$&$1$&$1$
\\ $35^*$&$0$&$0$&$0$&$0$&$0$&$1$&$0$&$0$&$0$&$0$&$1$&$0$&$0$&$1$&$0$&$1$
\\ $37$&$0$&$2$&$1$&$1$&$1$&$0$&$0$&$0$&$0$&$0$&$0$&$0$&$0$&$0$&$0$&$0$
\\ $71$&$1$&$1$&$1$&$1$&$1$&$0$&$0$&$0$&$1$&$0$&$0$&$0$&$0$&$0$&$0$&$0$
\\ $71^*$&$1$&$1$&$1$&$1$&$1$&$0$&$0$&$1$&$0$&$0$&$0$&$0$&$0$&$0$&$0$&$0$
\\ $117$&$0$&$3$&$1$&$1$&$2$&$0$&$0$&$0$&$0$&$0$&$0$&$0$&$0$&$0$&$0$&$0$
\\ $154$&$0$&$0$&$0$&$0$&$1$&$0$&$0$&$0$&$1$&$0$&$0$&$0$&$0$&$0$&$0$&$0$
\\ $154^*$&$0$&$0$&$0$&$0$&$1$&$0$&$0$&$1$&$0$&$0$&$0$&$0$&$0$&$0$&$0$&$0$
\\ $215$&$1$&$0$&$1$&$1$&$0$&$0$&$0$&$1$&$1$&$0$&$0$&$0$&$0$&$0$&$0$&$0$
\end{tabular}\end{center}

\begin{center}\begin{tabular}{c|cccccccccccccccc}
Block 3&$3$&$6$&$15_1$&$15_2$&$21$&$24$&$33$&$36$&$42$&$63$&$75$&$105$&$114$&$162$&$210$&$273$
\\\hline $3$&$0$&$0$&$0$&$0$&$0$&$0$&$1$&$1$&$0$&$0$&$1$&$0$&$1$&$1$&$0$&$0$
\\ $6$&$0$&$0$&$0$&$0$&$0$&$0$&$0$&$1$&$0$&$0$&$2$&$0$&$1$&$1$&$0$&$0$
\\ $15_1$&$0$&$0$&$0$&$0$&$0$&$0$&$1$&$1$&$0$&$1$&$0$&$0$&$0$&$1$&$0$&$1$
\\ $15_2$&$0$&$0$&$0$&$0$&$0$&$0$&$0$&$1$&$0$&$0$&$3$&$0$&$2$&$1$&$0$&$0$
\\ $21$&$0$&$0$&$0$&$0$&$0$&$0$&$1$&$0$&$0$&$0$&$0$&$0$&$0$&$0$&$1$&$0$
\\ $24$&$0$&$0$&$0$&$0$&$0$&$0$&$1$&$1$&$0$&$0$&$1$&$1$&$1$&$1$&$1$&$0$
\\ $33$&$1$&$0$&$1$&$0$&$1$&$1$&$0$&$0$&$0$&$0$&$0$&$0$&$0$&$0$&$0$&$0$
\\ $36$&$1$&$1$&$1$&$1$&$0$&$1$&$0$&$0$&$0$&$0$&$0$&$0$&$0$&$0$&$0$&$0$
\\ $42$&$0$&$0$&$0$&$0$&$0$&$0$&$0$&$0$&$0$&$1$&$0$&$1$&$1$&$1$&$0$&$0$
\\ $63$&$0$&$0$&$1$&$0$&$0$&$0$&$0$&$0$&$1$&$0$&$0$&$0$&$0$&$0$&$0$&$0$
\\ $75$&$1$&$2$&$0$&$3$&$0$&$1$&$0$&$0$&$0$&$0$&$0$&$0$&$0$&$0$&$0$&$0$
\\$105$&$0$&$0$&$0$&$0$&$0$&$1$&$0$&$0$&$1$&$0$&$0$&$0$&$0$&$0$&$0$&$0$
\\$114$&$1$&$1$&$0$&$2$&$0$&$1$&$0$&$0$&$1$&$0$&$0$&$0$&$0$&$0$&$0$&$0$
\\$162$&$1$&$1$&$1$&$1$&$0$&$1$&$0$&$0$&$1$&$0$&$0$&$0$&$0$&$0$&$0$&$0$
\\$210$&$0$&$0$&$0$&$0$&$1$&$1$&$0$&$0$&$0$&$0$&$0$&$0$&$0$&$0$&$0$&$0$
\\$273$&$0$&$0$&$1$&$0$&$0$&$0$&$0$&$0$&$0$&$0$&$0$&$0$&$0$&$0$&$0$&$0$
\end{tabular}\end{center}

\newpage
\section{\texorpdfstring{$\SU_3(3)$}{SU(3,3)}}

\subsection*{Block information}

\begin{center}\begin{tabular}{ccc}
\hline Block & Kernel & Modules
\\\hline Block 1 & $1$ & $1$, $3$, $3^*$, $6$, $6^*$, $7$, $15$, $15^*$
\\Block 2 & $1$ & $27$
\\\hline\end{tabular}\end{center}

\subsection*{Identification of modules}

\begin{enumerate}
\item Fix $3$ and $3^*$.
\item $6=S^2(3^*)$.
\item $15=\Lambda^2(6^*)$.
\end{enumerate}

\subsection*{Extension information}

\begin{center}\begin{tabular}{c|cccccccc}
Block 1&$1$&$3$&$3^*$&$6$&$6^*$&$7$&$15$&$15^*$
\\\hline $1$&$0$&$0$&$0$&$1$&$1$&$1$&$0$&$0$
\\ $3$&$0$&$0$&$1$&$0$&$0$&$1$&$0$&$1$
\\ $3^*$&$0$&$1$&$0$&$0$&$0$&$1$&$1$&$0$
\\ $6$&$1$&$0$&$0$&$0$&$0$&$0$&$0$&$1$
\\ $6^*$&$1$&$0$&$0$&$0$&$0$&$0$&$1$&$0$
\\ $7$&$1$&$1$&$1$&$0$&$0$&$0$&$0$&$0$
\\ $15$&$0$&$0$&$1$&$0$&$1$&$0$&$0$&$0$
\\ $15^*$&$0$&$1$&$0$&$1$&$0$&$0$&$0$&$0$
\end{tabular}\end{center}

\newpage
\section{\texorpdfstring{$\SU_3(4)$}{SU(3,4)}}


\subsection*{Block information}

\begin{center}\begin{tabular}{ccc}
\hline Block & Kernel & Modules
\\ \hline 1 & 1 & $1$, $3_1$, $3_1^*$, $3_2$, $3_2^*$, $8_1$, $8_2$, $9_{1,2}$, $9_{1,2}^*$, $\bar 9_{1,2}$, $\bar 9_{1,2}^*$, $24_{1,2}$, $24_{1,2}^*$, $24_{2,1}$, $24_{2,1}^*$
\\ 2 & 1 & $64$
\\\hline\end{tabular}\end{center}

\subsection*{Identification of modules}

We write $3_2$ for the twist of $3_1$, so that that $3_1$ and $3_2$ restrict to the same module for $\SL_3(2)$. Write $8_1$ for the module such that $8_1\otimes 3_2$ is irreducible.
\begin{enumerate}
\item $9_{1,2}=3_1\otimes 3_2$
\item $\bar 9_{1,2}=3_1\otimes 3_2^*$
\item $24_{1,2}=3_1\otimes 8_2$ and $24_{2,1}=3_2\otimes 8_1$
\end{enumerate}

\subsection*{Extension information}

\begin{center}\begin{tabular}{c|ccccccccccccccc}
Block 1& $1$&$3_1$&$3_1^*$&$3_2$&$3_2^*$&$8_1$&$8_2$&$9_{1,2}$&$9_{1,2}^*$&$\bar 9_{1,2}$&$\bar 9_{1,2}^*$&$24_{1,2}$&$24_{1,2}^*$&$24_{2,1}$&$24_{2,1}^*$
\\\hline $1$&$0$&$0$&$0$&$0$&$0$&$0$&$0$&$1$&$1$&$1$&$1$&$0$&$0$&$0$&$0$
\\ $3_1$&$0$&$0$&$0$&$1$&$1$&$0$&$0$&$0$&$1$&$0$&$0$&$0$&$0$&$1$&$0$
\\$3_1^*$&$0$&$0$&$0$&$1$&$1$&$0$&$0$&$1$&$0$&$0$&$0$&$0$&$0$&$0$&$1$
\\ $3_2$&$0$&$1$&$1$&$0$&$0$&$0$&$0$&$0$&$0$&$1$&$0$&$0$&$1$&$0$&$0$
\\ $3_2^*$&$0$&$1$&$1$&$0$&$0$&$0$&$0$&$0$&$0$&$0$&$1$&$1$&$0$&$0$&$0$
\\ $8_1$&$0$&$0$&$0$&$0$&$0$&$0$&$0$&$0$&$0$&$1$&$1$&$0$&$0$&$0$&$0$
\\ $8_2$&$0$&$0$&$0$&$0$&$0$&$0$&$0$&$1$&$1$&$0$&$0$&$0$&$0$&$0$&$0$
\\ $9_{1,2}$&$1$&$0$&$1$&$0$&$0$&$0$&$1$&$0$&$0$&$0$&$0$&$0$&$0$&$0$&$0$
\\ $9_{1,2}^*$&$1$&$1$&$0$&$0$&$0$&$0$&$1$&$0$&$0$&$0$&$0$&$0$&$0$&$0$&$0$
\\ $\bar 9_{1,2}$&$1$&$0$&$0$&$1$&$0$&$1$&$0$&$0$&$0$&$0$&$0$&$0$&$0$&$0$&$0$
\\ $\bar 9_{1,2}^*$&$1$&$0$&$0$&$0$&$1$&$1$&$0$&$0$&$0$&$0$&$0$&$0$&$0$&$0$&$0$
\\ $24_{1,2}$&$0$&$0$&$0$&$0$&$1$&$0$&$0$&$0$&$0$&$0$&$0$&$0$&$0$&$0$&$0$
\\ $24_{1,2}^*$&$0$&$0$&$0$&$1$&$0$&$0$&$0$&$0$&$0$&$0$&$0$&$0$&$0$&$0$&$0$
\\ $24_{2,1}$&$0$&$1$&$0$&$0$&$0$&$0$&$0$&$0$&$0$&$0$&$0$&$0$&$0$&$0$&$0$
\\ $24_{2,1}^*$&$0$&$0$&$1$&$0$&$0$&$0$&$0$&$0$&$0$&$0$&$0$&$0$&$0$&$0$&$0$
\end{tabular}\end{center}

\newpage
\section{\texorpdfstring{$\SU_3(5)$}{SU(3,5)}}

\subsection*{Block information}

\begin{center}\begin{tabular}{ccc}
\hline Block & Kernel & Modules
\\\hline 1 & 3 & $1$, $8$, $10$, $10^*$, $19$, $35$, $35^*$, $63$
\\ 2 & 3 & $125$
\\ 3 & 1 & $3$, $6$, $15_1$, $15_2$, $18$, $39$, $60$, $90$
\\ 4 & 1 & Dual to block 3
\\\hline\end{tabular}\end{center}

\subsection*{Identification of modules}

\begin{enumerate}
\item Fix $3$ and $3^*$. This fixes all modules in blocks 3 and 4 except for $15_1$ and $15_2$.
\item $10$ is a submodule of $3^*\otimes 6$.
\item $15_1$ is a submodule of $S^2(6^*)$.
\item $15_2=\Lambda^2(6^*)$.
\item $35$ is a submodule of $\Lambda^2(10^*)$.
\end{enumerate}

\subsection*{Extension information}

\begin{center}\begin{tabular}{c|cccccccc}
Block 1&$1$&$8$&$10$&$10^*$&$19$&$35$&$35^*$&$63$
\\ \hline $1$&$0$&$0$&$0$&$0$&$2$&$0$&$0$&$1$
\\ $8$&$0$&$0$&$0$&$0$&$3$&$1$&$1$&$2$
\\ $10$&$0$&$0$&$0$&$1$&$0$&$1$&$0$&$1$
\\ $10^*$&$0$&$0$&$1$&$0$&$0$&$0$&$1$&$1$
\\ $19$&$2$&$3$&$0$&$0$&$0$&$0$&$0$&$0$
\\ $35$&$0$&$1$&$1$&$0$&$0$&$0$&$0$&$0$
\\ $35^*$&$0$&$1$&$0$&$1$&$0$&$0$&$0$&$0$
\\ $63$&$1$&$2$&$1$&$1$&$0$&$0$&$0$&$0$
\end{tabular}\end{center}

\begin{center}\begin{tabular}{c|cccccccc}
Block 3&$3$&$6$&$15_1$&$15_2$&$18$&$39$&$60$&$90$
\\ \hline $3$&$0$&$0$&$0$&$0$&$2$&$3$&$0$&$0$
\\ $6$&$0$&$0$&$1$&$0$&$1$&$2$&$1$&$0$
\\ $15_1$&$0$&$1$&$0$&$0$&$0$&$0$&$0$&$1$
\\ $15_2$&$0$&$0$&$0$&$0$&$1$&$1$&$1$&$1$
\\ $18$&$2$&$1$&$0$&$1$&$0$&$0$&$0$&$0$
\\ $39$&$3$&$2$&$0$&$1$&$0$&$0$&$0$&$0$
\\ $60$&$0$&$1$&$0$&$1$&$0$&$0$&$0$&$0$
\\ $90$&$0$&$0$&$1$&$1$&$0$&$0$&$0$&$0$
\end{tabular}\end{center}

\newpage
\section{\texorpdfstring{$\SU_3(7)$}{SU(3,7)}}

\subsection*{Block information}

\begin{center}\begin{tabular}{ccc}
\hline Block & Kernel & Modules
\\\hline 1 & 1 & $1$, $3$, $3^*$, $6$, $6^*$, $8$, $10$, $10^*$, $15_1$, $15_1^*$, $15_2$, $15_2^*$, $21$, $21^*$, $24$, $24^*$, $27$, $28$, $28^*$,
\\ &&$33$, $33^*$, $35$, $35^*$, $36$, $36^*$, $37$, $42$, $42^*$, $63$, $63^*$, $71$, $71^*$, $75$, $75^*$, $105$, $105^*$,
\\ && $114$, $114^*$, $117$, $154$, $154^*$, $162$, $162^*$, $210$, $210^*$, $215$, $273$, $273^*$
\\ 2 & 1 & $343$
\\\hline\end{tabular}\end{center}

\subsection*{Identification of modules}

We produce an identification here that is consistent with the choice of highest weight labels from $\SL_3(7)$.

\begin{enumerate}
\item Fix $3$, $3^*$.
\item $6=S^2(3^*)$.
\item $10\mid 3^*\otimes 6$.
\item $15_1\mid S^2(6^*)$.
\item $15_2=\Lambda^2(6^*)$.
\item $21\oplus 24=3^*\otimes 15_1^*$.
\item $28\mid S^2(10^*)$.
\item $33$ is a factor of $3^*\otimes 21^*$.
\item $35\mid \Lambda^2(10^*)$.
\item $36 \subset S^2(21^*)$.
\item $42\oplus 63\mid \Lambda^2(15_1^*)$.
\item $71$ is a submodule of $8\otimes 28^*$.
\item $75$ is a factor of $6^*\otimes 42^*$.
\item $105\mid S^2(21^*)$.
\item $114$ lies inside $24\otimes 27$.
\item $154$ lies inside $\Lambda^2(28^*)$.
\item $162$ is a composition factor of $21^*\otimes 15_2^*$.
\item $210$ is a submodule of $15_1\otimes 28^*$.
\item $273$ is a submodule of $21\otimes 28$.
\end{enumerate}

\subsection*{Extension information}

\begin{center}
\begin{tabular}{c|c}
$M$&$\Ext^1=k$
\\\hline $1$& $36$, $36^*$, $215$.
\\$3$&$33^*$, $37$, $71^*$, $75^*$, $162$.
\\$6$&$36^*$, $71^*$, $75^*$, $114$, $117$.
\\$8$&$36$, $36^*$, $75$, $75^*$, $114$, $114^*$, $117$.
\\$10$&$33$, $71$, $75$, $114^*$, $162^*$.
\\$15_1$&$28^*$, $33$, $71$, $162^*$, $210^*$.
\\$15_2$&$37$, $71$, $75$, $75^*$, $114^*$, $117$, $162^*$, 
\\$21$&$21^*$, $63^*$, $215$.
\\$24$&$36$, $36^*$, $63^*$, $114^*$, $117$, $154$, $215$.
\\$27$&$37$, $75$, $75^*$, $105$, $105^*$, $162$, $162^*$.
\\$28$&$15_1^*$, $273$.
\\$33$&$3^*$, $10$, $15_1$, $35^*$.
\\$35$&$33^*$, $105^*$, $162$, $273$.
\\$36$&$1$, $6^*$, $8$, $24$, $24^*$.
\\$37$&$3$, $3^*$, $15_2$, $15_2^*$, $27$.
\\$42$& $71$, $114^*$, $154$, $210^*$.
\\$63$& $21^*$, $24^*$.
\\$71$& $3^*$, $6^*$, $10$, $15_1$, $15_2$, $42$.
\\$75$& $3^*$, $6^*$, $8$, $10$, $15_2$, $15_2^*$, $27$.
\\$105$& $27$, $35^*$.
\\$114$& $6$, $8$, $10^*$, $15_2^*$, $24^*$, $42^*$.
\\$117$& $6$, $6^*$, $8$, $15_2$, $15_2^*$, $24$, $24^*$.
\\$154$& $24$, $42$.
\\$162$& $3$, $10^*$, $15_1^*$, $15_2^*$, $27$, $35$.
\\$210$& $15_1^*$, $42^*$.
\\$215$& $1$, $21$, $21^*$, $24$, $24^*$.
\\$273$& $28$, $35$.
\end{tabular}\end{center}

\subsection*{Tensor products with the Steinberg}

\begin{center}
\begin{tabular}{c|c}
$M$&$M\otimes 343$
\\ \hline $3$&$P(273)$
\\$6$&$P(210)\oplus P(273)$
\\$8$&$P(215)\oplus 343^{\oplus 2}$
\\$10$&$P(154)\oplus P(215)\oplus 343$
\\$15_1$&$P(105)\oplus P(162)\oplus P(210)\oplus P(273)$
\\$15_2$&$P(162)\oplus P(210)\oplus P(273)^{\oplus 2}$
\\$21$&$P(63)\oplus P(114)\oplus P(162)\oplus P(210)\oplus P(273)$
\\$24$&$P(114)\oplus P(162)\oplus P(210)^{\oplus 2}\oplus P(273)^{\oplus 2}$
\\$27$&$P(117)\oplus P(154)\oplus P(154)^*\oplus P(215)^{\oplus 2}\oplus 343^{\oplus 3}$
\\$28$&$P(28)\oplus P(71)\oplus P(117)\oplus P(154)\oplus P(154^*)\oplus P(215)\oplus 343$
\\$33$&$P(33)\oplus P(75)\oplus P(105)\oplus P(114)\oplus P(162)\oplus P(210)\oplus P(273)$
\\$35$&$P(71^*)\oplus P(117)\oplus P(154)\oplus P(154^*)^{\oplus 2}\oplus P(215)^{\oplus 2}\oplus 343^{\oplus 2}$

\end{tabular}\end{center}

\newpage
\section{\texorpdfstring{$\SL_4(2)$}{SL(4,2)}}

\subsection*{Block information}

\begin{center}\begin{tabular}{ccc}
\hline Block & Kernel & Modules
\\\hline Block 1 & $1$ & $1$, $4$, $4^*$, $6$, $14$, $20$, $20^*$
\\Block 2 & $1$ & $64$
\\\hline\end{tabular}\end{center}

\subsection*{Identification of modules}
\begin{enumerate}
\item Fix $4$.
\item $20$ is a submodule of $4\otimes 6$.
\end{enumerate}

\subsection*{Highest weight labels}

\begin{center}\begin{tabular}{ccccccccc}
\hline $\lambda$&000&001&100&010&101&011&110&111
\\\hline $L(\lambda)$ &$1$&$4$&$4^*$&$6$&$14$&$20$&$20^*$&$64$
\\ $\dim(W(\lambda))$&$1$&$4$&$4$&$6$&$15$&$20$&$20$&$64$
\\\hline
\end{tabular}\end{center}

\subsection*{Extension information}

\begin{center}\begin{tabular}{c|ccccccc}
Block 1&$1$&$4$&$4^*$&$6$&$14$&$20$&$20^*$
\\\hline $1$&$0$&$0$&$0$&$1$&$1$&$1$&$1$
\\ $4$&$0$&$0$&$1$&$1$&$0$&$1$&$0$
\\ $4$&$0$&$1$&$0$&$1$&$0$&$0$&$1$
\\ $6$&$1$&$1$&$1$&$0$&$1$&$0$&$0$
\\ $14$&$1$&$0$&$0$&$1$&$0$&$0$&$0$
\\ $20$&$1$&$1$&$0$&$0$&$0$&$0$&$0$
\\ $20^*$&$1$&$0$&$1$&$0$&$0$&$0$&$0$
\end{tabular}\end{center}

\newpage
\section{\texorpdfstring{$\SL_4(3)$}{SL(4,3)}}

\subsection*{Block information}

\begin{center}\begin{tabular}{ccc}
\hline Block & Kernel & Modules
\\\hline Block 1 & $2$ & $1$, $6$, $10$, $10^*$, $15$, $19$, $44$, $45$, $45^*$, $69$, $126$, $126^*$, $156$, $294$
\\Block 2 & $2$ & $729$
\\ Block 3& $1$ & $4$, $4^*$, $16$, $16^*$, $36$, $36^*$, $60$, $60^*$, $116$, $116^*$, $360$, $360^*$
\\\hline\end{tabular}\end{center}

\subsection*{Identification of modules}
\begin{enumerate}
\item Fix $4$.
\item $S^2(4)=10$.
\item $S^3(4)=16/4$.
\item $4\otimes 15=4/16^*/4\oplus 36$.
\item $4\otimes 19=16^*\oplus 60^*$.
\item $16\otimes 10^*=36\oplus 4/116^*/4$.
\item $4\otimes 729=P(360)$.
\item $\Lambda^2(10)=45^*$.
\item $126^*\mid S^3(10)$.
\end{enumerate}

\subsection*{Highest weight labels}

\begin{center}\begin{tabular}{cccccccccccccccc}
\hline $\lambda$&000&010&200&002&101&020&111&012&210&202&022&220&121&212&222
\\\hline $L(\lambda)$ &$1$&$6$&$10$&$10^*$&$15$&$19$&$44$&$45$&$45^*$&$69$&$126$&$126^*$&$156$&$294$&$729$
\\ $\dim(W(\lambda))$&$1$&$6$&$10$&$10$&$15$&$20$&$64$&$45$&$45$&$84$&$126$&$126$&$175$&$300$&$729$
\\\hline
\end{tabular}\end{center}

\begin{center}\begin{tabular}{ccccccccccccc}
\hline $\lambda$&100&001&110&011&201&102&021&120&211&112&221&122
\\\hline $L(\lambda)$ &$4$&$4^*$&$16$&$16^*$&$36$&$36^*$&$60$&$60^*$&$116$&$116^*$&$360$&$360^*$
\\ $\dim(W(\lambda))$&$4$&$4$&$20$&$20$&$36$&$36$&$60$&$60$&$140$&$140$&$360$&$360$
\\\hline
\end{tabular}\end{center}

\subsection*{Extension information}

\begin{center}\begin{tabular}{c|cccccccccccccc}
Block 1&$1$&$6$&$10$&$10^*$&$15$&$19$&$44$&$45$&$45^*$&$69$&$126$&$126^*$&$156$&$294$
\\\hline $1$&$0$&$0$&$0$&$0$&$0$&$1$&$2$&$0$&$0$&$1$&$0$&$0$&$0$&$0$
\\$6$&$0$&$0$&$0$&$0$&$0$&$2$&$1$&$0$&$0$&$0$&$0$&$0$&$0$&$2$
\\$10$&$0$&$0$&$0$&$0$&$0$&$1$&$1$&$1$&$0$&$0$&$0$&$0$&$1$&$1$
\\$10^*$&$0$&$0$&$0$&$0$&$0$&$1$&$1$&$0$&$1$&$0$&$0$&$0$&$1$&$1$
\\$15$&$0$&$0$&$0$&$0$&$0$&$0$&$2$&$0$&$0$&$2$&$1$&$1$&$0$&$0$
\\$19$&$1$&$2$&$1$&$1$&$0$&$0$&$0$&$0$&$0$&$0$&$0$&$0$&$1$&$0$
\\$44$&$2$&$1$&$1$&$1$&$2$&$0$&$0$&$0$&$0$&$0$&$0$&$0$&$0$&$0$
\\$45$&$0$&$0$&$1$&$0$&$0$&$0$&$0$&$0$&$0$&$0$&$1$&$0$&$0$&$0$
\\$45^*$&$0$&$0$&$0$&$1$&$0$&$0$&$0$&$0$&$0$&$0$&$0$&$1$&$0$&$0$
\\$69$&$1$&$0$&$0$&$0$&$2$&$0$&$0$&$0$&$0$&$0$&$0$&$0$&$0$&$0$
\\$126$&$0$&$0$&$0$&$0$&$1$&$0$&$0$&$1$&$0$&$0$&$0$&$0$&$0$&$0$
\\$126^*$&$0$&$0$&$0$&$0$&$1$&$0$&$0$&$0$&$1$&$0$&$0$&$0$&$0$&$0$
\\$156$&$0$&$0$&$1$&$1$&$0$&$1$&$0$&$0$&$0$&$0$&$0$&$0$&$0$&$0$
\\$294$&$0$&$2$&$1$&$1$&$0$&$0$&$0$&$0$&$0$&$0$&$0$&$0$&$0$&$0$
\end{tabular}\end{center}

\begin{center}\begin{tabular}{c|cccccccccccc}
Block 3&$4$&$4^*$&$16$&$16^*$&$36$&$36^*$&$60$&$60^*$&$116$&$116^*$&$360$&$360^*$
\\ \hline $4$&$0$&$0$&$1$&$1$&$0$&$0$&$2$&$0$&$2$&$2$&$0$&$0$
\\ $4^*$&$0$&$0$&$1$&$1$&$0$&$0$&$0$&$2$&$2$&$2$&$0$&$0$
\\ $16$&$1$&$1$&$0$&$1$&$0$&$0$&$0$&$0$&$0$&$0$&$1$&$0$
\\ $16^*$&$1$&$1$&$1$&$0$&$0$&$0$&$0$&$0$&$0$&$0$&$0$&$1$
\\ $36$&$0$&$0$&$0$&$0$&$0$&$1$&$1$&$0$&$1$&$0$&$0$&$0$
\\ $36^*$&$0$&$0$&$0$&$0$&$1$&$0$&$0$&$1$&$0$&$1$&$0$&$0$
\\ $60$&$2$&$0$&$0$&$0$&$1$&$0$&$0$&$0$&$0$&$0$&$0$&$0$
\\ $60^*$&$0$&$2$&$0$&$0$&$0$&$1$&$0$&$0$&$0$&$0$&$0$&$0$
\\ $116$&$2$&$2$&$0$&$0$&$1$&$0$&$0$&$0$&$0$&$0$&$0$&$0$
\\ $116^*$&$2$&$2$&$0$&$0$&$0$&$1$&$0$&$0$&$0$&$0$&$0$&$0$
\\ $360$&$0$&$0$&$1$&$0$&$0$&$0$&$0$&$0$&$0$&$0$&$0$&$0$
\\ $360^*$&$0$&$0$&$0$&$1$&$0$&$0$&$0$&$0$&$0$&$0$&$0$&$0$
\end{tabular}\end{center}

\newpage
\section{\texorpdfstring{$\SL_4(5)$}{SL(4,5)}}

\subsection*{Block information}

\begin{center}\begin{tabular}{ccc}
\hline Block & Kernel & Modules
\\\hline Block 1 & 4 & $1$, $15$, $20_1$, $35$, $35^*$, $45$, $45^*$, $83$, $85$, $173$, $173^*$, $175$, $235$, $235^*$, $299$, $300_1$, $360$, $360^*$, $525$,
\\ && $720_1$, $875$, $875^*$, $885$, $885^*$, $1375$, $1375^*$, $1396$, $1920$, $3420$, $3420^*$, $5025$, $6935$
\\Block 2 & 4 & $15625$
\\ Block 3& 2 & $6$, $10$, $10^*$, $50$, $58$, $68$, $68^*$, $70$, $70^*$, $140_1$, $140_1^*$, $184$, $270$, $270^*$, $300_2$, $420_1$, $530$, $530^*$,
\\ && $560$, $560^*$, $750$, $750^*$, $1190$, $1270$, $1270^*$, $1780$, $1780^*$, $2150$, $3686$, $6500$, $6500^*$, $9400$
\\ Block 4 & 1 &  $4$, $20_2$, $20_3$, $32$, $52$, $60$, $80$, $120$, $140_2$, $160$, $192$, $220$, $260$, $380$, $420_2$, $460$, $476$, $500$, $720_2$,
\\ && $900$, $960$, $1000$, $1400$, $1420$, $1684$, $2940$, $3500$, $4060$, $6060$, $10500$
\\ Block 5 & 1 & Dual to Block 4
\\\hline\end{tabular}\end{center}

\subsection*{Identification of modules}

Fix $4$.

\begin{enumerate}
\item $S^2(4)=10$.
\item $S^2(20_2)$ has $68$ as a submodule.
\item $\Lambda^2(20_3)$ has $70$ as a summand.
\item $\Lambda^2(20_2)$ has $140_1$ as a summand.
\end{enumerate}

\subsection*{Extension information}

\begin{center}
\begin{tabular}{ccccc}\hline
Module & Tested to & $\Ext^1=k$ & $\Ext^1=k^2$
\\ \hline $1$ & All & $83$, $173$, $173^*$ & $299$, $1396$
\\ $15$ & $360$ & $173$, $173^*$, $235$, $235^*$, $299$ & None
\\ $20_1$ & $300$& $85$, $235$, $235^*$, $299$ & None
\\ $35$ & $175$ & $85$ & None
\\ $45$ & $175$ & $85$ & None
\\ $83$ & $85$ & $1$ & None
\\ $85$ & $85$ & $20_1$, $35$, $35^*$, $45$, $45^*$ & None
\\ \hline\end{tabular}
\end{center}

\begin{center}
\begin{tabular}{cccc}\hline
Module & Tested to & $\Ext^1=k$ & $\Ext^1=k^2$
\\ \hline $6$ & $184$ & $58$, $68$, $68^*$ & $184$
\\ $10$ & $184$ & $68^*$ & None
\\ $50$ & $184$ & None & None
\\ $58$ & $140$ & $6$, $68$, $68^*$  & None
\\ $68$ & $70$ & $6$, $10^*$, $58$, $68^*$ & None 
\\ $70$ & $70$ & None & None 
\\ \hline\end{tabular}
\end{center}

\begin{center}
\begin{tabular}{cccc}\hline
Module & Tested to & $\Ext^1=k$ & $\Ext^1=k^2$
\\ \hline $4$ & $1684$ & $32$, $52$, $192$, $220$, $380$, $720_2$, $960$ & $476$, $1684$
\\ $20_2$ & $160$ & $80$ & None
\\ $20_3$ & $160$ & $80$ & None
\\ $32$ & $140$ & $4$, $52$ & None
\\ $52$ & $140$ & $4$, $32$, $80$ & None
\\ \hline\end{tabular}
\end{center}

\newpage
\section{\texorpdfstring{$\SL_4(7)$}{SL(4,7)}}

\subsection*{Block information}

These modules are only necessarily complete up to dimension 100.

\begin{center}\begin{tabular}{ccc}
\hline Block & Kernel & Modules
\\\hline1 & 2 & $1$, $6$, $10$, $10^*$, $15$, $20_1$, $35$, $35^*$, $45$, $45^*$, $50$, $64$, $70$, $70^*$, $84_1$, $84_2$, $84_2^*$, $105$, $126$, $126^*$,
\\&& $140_1$, $140_1^*$, $149$, $149^*$, $175$, $189$, $189^*$, $196$, $206$, $206^*$, $211$, $211^*$, $231$, $236$, $260$, $260^*$,
\\&& $280_1$, $280_1^*$, $285$, $315$, $315^*$, $334$, $380$, $380^*$, $484$, $514$, $514^*$, $571$, $571^*$, $616$, $616^*$, $630$, $630^*$,
\\&& $665$, $665^*$, $709$, $735$, $769$, $770_1$, $770_1^*$, $770_2$, $770_2^*$, $824$, $875$, $875^*$, $954$, $1026$, $1026^*$, $1084$,
\\&& $1164$, $1164^*$, $1200$, $1200^*$, $1414$, $1414^*$, $1470$, $1480$, $1480^*$, $1540$, $1540^*$, $1680$, $1680^*$,
\\&& $1820$, $1820^*$, \dots
\\2 & 2 & $117649$
\\3 & 1 & $4$, $4^*$, $20_2$, $20_2^*$, $20_3$, $20_3^*$, \dots
\\\hline\end{tabular}\end{center}

\subsection*{Identification of modules}

\begin{enumerate}
\item Fix $4$ and $4^*$.
\item $S^2(4)=10$.
\item $S^4(4)=35$.
\item $149$ is a composition factor of $20_1\otimes 35$.
\end{enumerate}

\subsection*{Extension information}

\begin{center}\begin{tabular}{ccc}\hline
Module & Tested to & $\Ext^1=k$
\\\hline $1$&$571$&$380$, $380^*$ 
\\$6$&$236$&$149$, $149^*$
\\$10$&$236$&$149$
\\$15$&$236$&None
\\$20_1$&$236$&None
\\$35$&$196$&None
\\$45$&$84$&None
\\$50$&$84$&None
\\$64$&$84$&None
\\$70$&$70$&None
\\$84_1$ & $64$ & None
\\$84_2$ & $64$ & None
\\ \hline\end{tabular}\end{center}

\newpage
\section{\texorpdfstring{$\SU_4(2)$}{SU(4,2)}}

\subsection*{Block information}

\begin{center}\begin{tabular}{ccc}
\hline Block & Kernel & Modules
\\\hline $1$ & 1 & $1$, $4$, $4^*$, $6$, $14$, $20$, $20^*$
\\ $2$ & 1 & $64$
\\\hline\end{tabular}\end{center}

\subsection*{Extension information}

\begin{center}\begin{tabular}{c|ccccccc}
Block 1&$1$&$4$&$4^*$&$6$&$14$&$20$&$20^*$
\\\hline $1$&$0$&$1$&$1$&$0$&$1$&$0$&$0$
\\$4$&$1$&$0$&$0$&$1$&$0$&$0$&$0$
\\$4^*$&$1$&$0$&$0$&$1$&$0$&$0$&$0$
\\$6$&$0$&$1$&$1$&$0$&$1$&$1$&$1$
\\$14$&$1$&$0$&$0$&$1$&$0$&$0$&$0$
\\$20$&$0$&$0$&$0$&$1$&$0$&$0$&$0$
\\$20^*$&$0$&$0$&$0$&$1$&$0$&$0$&$0$
\end{tabular}\end{center}

\subsection*{Actions of unipotent elements on modules}

\begin{center}\begin{tabular}{c|ccccccc}
Block 1&$1$&$4$&$4^*$&$6$&$14$&$20$&$20^*$
\\\hline 2A&$1$&$2,1^2$&$2,1^2$&$2^2,1^2$&$2^6,1^2$&$2^9,1^2$&$2^9,1^2$
\\ 2B&$1$&$2^2$&$2^2$&$2^2,1^2$&$2^6,1^2$&$2^{10}$&$2^{10}$
\\ 4A[A]&$1$&$3,1$&$3,1$&$3^2$&$4^2,3^2$&$4^4,3,1$&$4^4,3,1$
\\ 4B[B]&$1$&$4$&$4$&$4,2$&$4^3,2$&$4^5$&$4^5$
\end{tabular}\end{center}

\newpage
\section{\texorpdfstring{$\SU_4(3)$}{SU(4,3)}}

\subsection*{Block information}

\begin{center}\begin{tabular}{ccc}
\hline Block & Kernel & Modules
\\\hline 1 & $4$ & $1$, $15$, $19$, $45$, $45^*$, $69$, $156$
\\ 2 & $4$ & $729$
\\ 3 & $2$ & $6$, $10$, $10^*$, $44$, $126$, $126^*$, $294$
\\ 4 & $1$ & $4$, $16$, $36$, $60$, $116$, $360$
\\ 5 & $1$ & Dual to block 4
\\\hline\end{tabular}\end{center}

\subsection*{Identification of modules}

\begin{enumerate}
\item $S^2(16)=10\oplus 126$.
\item $\Lambda^2(10)=45$.
\end{enumerate}
\subsection*{Extension information}

\begin{center}\begin{tabular}{c|ccccccc}
Block 1&$1$&$15$&$19$&$45$&$45^*$&$69$&$156$
\\\hline $1$&$0$&$0$&$3$&$0$&$0$&$1$&$0$
\\$15$&$0$&$0$&$2$&$1$&$1$&$2$&$0$
\\$19$&$3$&$2$&$0$&$0$&$0$&$0$&$3$
\\ $45$&$0$&$1$&$0$&$0$&$1$&$0$&$1$
\\ $45^*$&$0$&$1$&$0$&$1$&$0$&$0$&$1$
\\$69$&$1$&$2$&$0$&$0$&$0$&$0$&$0$
\\$156$&$0$&$0$&$3$&$1$&$1$&$0$&$0$
\end{tabular}\end{center}

\begin{center}\begin{tabular}{c|ccccccc}
Block 3&$6$&$10$&$10^*$&$44$&$126$&$126^*$&$294$
\\\hline $6$&$0$&$0$&$0$&$3$&$0$&$0$&$2$
\\ $10$&$0$&$0$&$1$&$2$&$0$&$1$&$1$
\\ $10^*$&$0$&$1$&$0$&$2$&$1$&$0$&$1$
\\ $44$&$3$&$2$&$2$&$0$&$0$&$0$&$0$
\\ $126$&$0$&$0$&$1$&$0$&$0$&$0$&$0$
\\ $126^*$&$0$&$1$&$0$&$0$&$0$&$0$&$0$
\\ $294$&$2$&$1$&$1$&$0$&$0$&$0$&$0$
\end{tabular}\end{center}

\begin{center}\begin{tabular}{c|cccccc}
Block 4&$4$&$16$&$36$&$60$&$116$&$360$
\\\hline $4$&$0$&$3$&$0$&$2$&$3$&$0$
\\ $16$&$3$&$0$&$1$&$0$&$1$&$0$
\\ $36$&$0$&$1$&$0$&$1$&$0$&$1$
\\ $60$&$2$&$0$&$1$&$0$&$0$&$1$
\\ $116$&$3$&$1$&$0$&$0$&$0$&$0$
\\ $360$&$0$&$0$&$1$&$1$&$0$&$0$
\end{tabular}\end{center}

\subsection*{Tensor products with the Steinberg}

\begin{center}\begin{tabular}{c|c}
$M$ & $M\otimes 729$
\\\hline $15$ & $P(156)\oplus 729^{\oplus 3}$
\\ $19$ & $P(69)\oplus P(156)\oplus 729^{\oplus 2}$
\\ \hline $6$ & $P(294)$
\\ $10$ & $P(126^*)\oplus P(294)$
\\ \hline $4$ & $P(360)$
\\ $16$ & $P(116)\oplus P(360)$
\end{tabular}\end{center}

\newpage
\section{\texorpdfstring{$\SU_4(5)$}{SU(4,5)}}

\subsection*{Block information}
\begin{center}\begin{tabular}{ccc}
\hline Block & Kernel & Modules
\\\hline Block 1 & $2$ & $1,6,10,10^*,15,20_1,35,35^*,45,45^*,50,58,68,68^*,70,70^*,83,85,140_1,140_1^*,173,173^*,$
\\ &&$175,184,235,235^*,270,270^*,299,300_1,300_2,360,360^*,420_1,525,530,530^*$,
\\ && $560,560^*,720_1,750,750^*,875,875^*,885,885^*,1190,1270,1270^*,1375,1375^*,1396,$
\\&&$1780,1780^*,1920,2150,3420,3420^*,3686,5025,6500,6500^*,6935,9400$
\\Block 2 & $2$ & $15625$
\\ Block 3 & $1$ &  $4,4^*,20_2,20_2^*,20_3,20_3^*,32,32^*,52,52^*,60,60^*,80,80^*,120,120^*,140_2,140_2^*,160,160^*,$
\\ &&$192,192^*,220,220^*,260,260^*,380,380^*,420_2,420_2^*,460,460^*,476,476^*,500,500^*,$
\\ && $720_2,720_2^*,900,900^*,960,960^*,1000,1000^*,1400,1400^*,1420,1420^*,1684,1684^*,$
\\ &&$2940,2940^*,3500,3500^*,4060,4060^*,6060,6060^*,10500,10500^*$
\\\hline\end{tabular}\end{center}

\subsection*{Extension information}

\begin{center}\begin{tabular}{cccc}\hline
Module & Tested to & $\Ext^1=k$ & $\Ext^1=k^2$
\\\hline $1$ & $300$ & $68,68^*,83$ & $299$
\\ $6$ & $68$ & $58$ & None
\\ $10$ & $68$ & None & None
\\ $15$ & $68$ & $68,68^*$ & None
\\$20_1$ & $68$ & None & None
\\ $35$ & $58$ & None & None
\\ $45$ & $58$ & None & None
\\ $50$ & $58$ & None & None
\\ $58$ & $58$ & $6$ & None
\\ \hline\end{tabular}\end{center}

\newpage
\section{\texorpdfstring{$\SU_4(7)$}{SU(4,7)}}

\subsection*{Block information}
\begin{center}\begin{tabular}{ccc}
\hline Block & Kernel & Modules
\\\hline Block 1 & $4$ & $1,15,20_1,35,35^*,45,45^*,84_1,105,149,149^*,175,189,189^*,211,211^*,231$,
\\&&$280_1,280_1^*,285,315,315^*,571,571^*,616,616^*,665,665^*,709,735,824,875,875^*$,
\\&&$1164,1164^*,1200,1200^*,1480,1480^*,1540,1540^*,2156...$
\\Block 2 & $4$ & $117649$
\\ Block 3& $2$ & $6,10,10^*,50,64,70,70^*,84_2,84_2^*,126,126^*,140,140^*,196,206,206^*,236,260,260^*$,
\\&&$334,380,380^*,514,514^*,630,630^*,770_1,770_1^*,770_2,770_2^*,954...$
\\ Block 4 & $1$ & $4,20_2,20_3,36,56,60,84_3,100,116,124,140_1,140_2,180,224,280_2,280_3,300$,
\\&&$416,420_1,420_2,480,496,520,756,924,...$
\\ Block 5 & $1$ & Dual to Block 4
\\\hline\end{tabular}\end{center}

\subsection*{Extension information}

\begin{center}\begin{tabular}{ccc}
\hline Module & Tested to & $\Ext^1=k$ 
\\\hline $1$ & $571$ & $149,149^*$
\\$15$ & $84$ & None
\\$20$ & $45$ & None
\\$35$ & $45$ & None
\\$45$ & $45$ & None
\\\hline\end{tabular}\end{center}

\newpage
\section{\texorpdfstring{$\SL_5(2)$}{SL(5,2)}}

\subsection*{Block information}

\begin{center}\begin{tabular}{ccc}
\hline Block & Kernel & Modules
\\\hline Block 1 & $1$ & $1$, $5$, $5^*$, $10$, $10^*$, $24$, $40_1$, $40_1^*$, $40_2$, $40_2^*$, $74$, $160$, $160^*$, $280$, $280^*$
\\Block 2 & $1$ & $1024$
\\\hline\end{tabular}\end{center}

\subsection*{Identification of modules}

Fix $5$ and $5^*$.

\begin{enumerate}
\item $\Lambda^2(5^*)=10$;
\item $\Lambda^2(10)=5\oplus 40_1$;
\item $5\otimes 10^*=10\oplus 40_2$;
\item $160$ is a composition factor of $10\otimes 24$;
\item $280$ is a composition factor of $5\otimes 74$.
\end{enumerate}

\subsection*{Highest weight labels}

\begin{center}\begin{tabular}{ccccccccc}
\hline $\lambda$ &0000&1000&0001&0010&0100&1001&0101&1010
\\\hline $L(\lambda)$&$1$&$5$&$5^*$&$10$&$10^*$&$24$&$40_1$&$40_1^*$
\\ $\dim(W(\lambda))$&$1$&$5$&$5$&$10$&$10$&$24$&$45$&$45$
\\ \hline $\lambda$&1100&0011&0110&1101&1011&1110&0111&1111
\\\hline $L(\lambda)$ &$40_2$&$40_2^*$&$74$&$160$&$160^*$&$280$&$280^*$&$1024$
\\ $\dim(W(\lambda))$&$40$&$40$&$75$&$175$&$175$&$280$&$280$&$1024$
\\ \hline
\end{tabular}\end{center}

\subsection*{Extension information}

\begin{center}\begin{tabular}{c|ccccccccccccccc}
Block 1&$1$&$5$&$5^*$&$10$&$10^*$&$24$&$40_1$&$40_1^*$&$40_2$&$40_2^*$&$74$&$160$&$160^*$&$280$&$280^*$
\\\hline $1$&$0$&$0$&$0$&$0$&$0$&$0$&$1$&$1$&$0$&$0$&$1$&$0$&$0$&$0$&$0$
\\ $5$&$0$&$0$&$0$&$0$&$1$&$0$&$1$&$0$&$0$&$1$&$1$&$1$&$1$&$0$&$0$
\\ $5^*$&$0$&$0$&$0$&$1$&$0$&$0$&$0$&$1$&$1$&$0$&$1$&$1$&$1$&$0$&$0$
\\ $10$&$0$&$0$&$1$&$0$&$1$&$0$&$1$&$0$&$0$&$0$&$0$&$0$&$0$&$0$&$1$
\\ $10^*$&$0$&$1$&$0$&$1$&$0$&$0$&$0$&$1$&$0$&$0$&$0$&$0$&$0$&$1$&$0$
\\ $24$&$0$&$0$&$0$&$0$&$0$&$0$&$1$&$1$&$1$&$1$&$0$&$0$&$0$&$0$&$0$
\\ $40_1$&$1$&$1$&$0$&$1$&$0$&$1$&$0$&$0$&$0$&$0$&$0$&$0$&$0$&$0$&$0$
\\ $40_1^*$&$1$&$0$&$1$&$0$&$1$&$1$&$0$&$0$&$0$&$0$&$0$&$0$&$0$&$0$&$0$
\\ $40_2$&$0$&$0$&$1$&$0$&$0$&$1$&$0$&$0$&$0$&$0$&$0$&$0$&$0$&$1$&$0$
\\ $40_2^*$&$0$&$1$&$0$&$0$&$0$&$1$&$0$&$0$&$0$&$0$&$0$&$0$&$0$&$0$&$1$
\\ $74$&$1$&$1$&$1$&$0$&$0$&$0$&$0$&$0$&$0$&$0$&$0$&$0$&$0$&$0$&$0$
\\ $160$&$0$&$1$&$1$&$0$&$0$&$0$&$0$&$0$&$0$&$0$&$0$&$0$&$0$&$0$&$0$
\\ $160^*$&$0$&$1$&$1$&$0$&$0$&$0$&$0$&$0$&$0$&$0$&$0$&$0$&$0$&$0$&$0$
\\ $280$&$0$&$0$&$0$&$0$&$1$&$0$&$0$&$0$&$1$&$0$&$0$&$0$&$0$&$0$&$0$
\\ $280^*$&$0$&$0$&$0$&$1$&$0$&$0$&$0$&$0$&$0$&$1$&$0$&$0$&$0$&$0$&$0$
\end{tabular}\end{center}

\newpage
\section{\texorpdfstring{$\SL_5(3)$}{SL(5,3)}}

\subsection*{Block information}

\begin{center}\begin{tabular}{ccc}
\hline Block & Kernel & Modules
\\\hline $1$ & $1$ & $1,5,5^*,10,10^*,15,15^*,24,30,30^*,45_1,45_1^*, 45_2,45_2^*,51,65,65^*,105,105^*,$
\\&&$126,126^*,135,135^*,165,165^*,174,174^*,195,195^*,199,315,315^*,375,375^*,390,390^*,$
\\ && $420,420^*,450,450^*,476,510,510^*,720,720^*,835,835^*,855,855^*,1125,1215,1215^*,$
\\&&$1440,1440^*,1545,1545^*,1701,1701^*,2430,2430^*,2535,2535^*,2574,2574^*,3465,3465^*,$
\\&&$3780,3780^*,4199,6810,6810^*,7350,7350^*,8505,8505^*,10824,19260,19260^*,26730,26730^*$
\\ $2$ & $1$ & $59049$
\\\hline\end{tabular}\end{center}

\subsection*{Extension information}

\begin{center}\begin{tabular}{cccc}\hline
Module & Tested to & $\Ext^1=k$ & $\Ext^1=k^2$
\\\hline $1$ & $6810$ & $135,135^*,174,174^*,199,510,510^*,835,835^*,2535,2535^*$ & $4199$
\\ $5$ & $65$ & $30,45_1,65$ & None
\\ $10$ & $51$ & $30$, $45_1$ & None
\\ $15$ & $51$ & $45_1$ & None
\\ $24$ & $51$ & $51$ & None
\\ $30$ & $51$ & $5$, $10$, $51$ & None
\\ $45_1$ & $51$ & $5,10,15,45_1^*$ & None
\\ $45_2$ & $51$ & None & None
\\ $51$ & $51$ & $24,30,30^*$ & None
\\ \hline\end{tabular}\end{center}

\newpage
\section{\texorpdfstring{$\SL_5(5)$}{SL(5,5)}}

\subsection*{Block information}

These modules are only necessarily complete up to dimension 1200.

\begin{center}\begin{tabular}{ccc}
\hline Block & Kernel & Modules
\\\hline $1$ & $1$ & $1,5,5^*,10,10^*,15,15^*,23,35,35^*,40,40^*,45,45^*,50,50^*,70_1,70_1^*,70_2,70_2^*,75,$
\\&&$103,103^*,105,105^*,121,121^*,160,160^*,170,170^*,175_1,175_1^*,175_2,175_2^*,175_3,175_3^*,$
\\&&$185,185^*,200,235,235^*,255,255^*,280,280^*,305,305^*,315,315^*,320,320^*,365,365^*,381,$
\\&&$410,410^*,420,420^*,450_1,450_1^*,450_2,450_2^*,535,535^*,540,540^*,545,545^*,615,615^*,670,670^*,$
\\ &&$683,683^*,700_1,700_1^*,700_2,700_2^*,794,875,875^*,895,895^*,955,955^*,1000,1050_1,1050_1^*$
\\ &&$1050_2,1050_2^*,1055,1055^*,1063,1063^*,1090,1090^*,1115,1115^*,1150,1150^*,1200,1200^*$,...
\\ $2$ & $1$ & $9765625$
\\\hline\end{tabular}\end{center}

\subsection*{Extension information}

\begin{center}\begin{tabular}{ccc}\hline
Module & Tested to & $\Ext^1=k$ 
\\\hline $1$ & $1150$ & $23,535,535^*,794$
\\ $5$ & $121$ & $121$
\\ $10$ & $50$ & None
\\ $15$ & $50$ & None
\\ $23$ & $50$ & $1$
\\ $35$ & $50$ & None
\\ $40$ & $50$ & None
\\ $45$ & $50$ & None
\\ $50$ & $50$ & None
\\ \hline\end{tabular}\end{center}

\newpage
\section{\texorpdfstring{$\SU_5(2)$}{SU(5,2)}}

\subsection*{Block information}

\begin{center}\begin{tabular}{ccc}
\hline Block & Kernel & Modules
\\\hline Block 1 & $1$ & $1$, $5$, $5^*$, $10$, $10^*$, $24$, $40_1$, $40_1^*$, $40_2$, $40_2^*$, $74$, $160$, $160^*$, $280$, $280^*$
\\Block 2 & $1$ & $1024$
\\\hline\end{tabular}\end{center}

\subsection*{Identification of modules}

Fix $5$ and $5^*$.

\begin{enumerate}
\item $\Lambda^2(5^*)=10$;
\item $\Lambda^2(10)=5\oplus 40_1$;
\item $5\otimes 10^*=10\oplus 40_2$;
\item $160$ is a composition factor of $10\otimes 24$;
\item $280$ is a submodule of $5\otimes 74$.
\end{enumerate}

\subsection*{Extension information}

\begin{center}\begin{tabular}{c|ccccccccccccccc}
Block 1&$1$&$5$&$5^*$&$10$&$10^*$&$24$&$40_1$&$40_1^*$&$40_2$&$40_2^*$&$74$&$160$&$160^*$&$280$&$280^*$
\\\hline $1$&$0$&$0$&$0$&$1$&$1$&$0$&$0$&$0$&$0$&$0$&$1$&$0$&$0$&$0$&$0$
\\$5$&$0$&$0$&$1$&$1$&$0$&$0$&$1$&$0$&$1$&$0$&$0$&$1$&$0$&$0$&$0$
\\$5^*$&$0$&$1$&$0$&$0$&$1$&$0$&$0$&$1$&$0$&$1$&$0$&$0$&$1$&$0$&$0$
\\$10$&$1$&$1$&$0$&$0$&$0$&$1$&$1$&$1$&$0$&$0$&$0$&$0$&$1$&$0$&$0$
\\$10^*$&$1$&$0$&$1$&$0$&$0$&$1$&$1$&$1$&$0$&$0$&$0$&$1$&$0$&$0$&$0$
\\$24$&$0$&$0$&$0$&$1$&$1$&$0$&$0$&$0$&$0$&$0$&$0$&$0$&$0$&$0$&$0$
\\$40_1$&$0$&$1$&$0$&$1$&$1$&$0$&$0$&$0$&$0$&$0$&$0$&$0$&$0$&$0$&$1$
\\$40_1^*$&$0$&$0$&$1$&$1$&$1$&$0$&$0$&$0$&$0$&$0$&$0$&$0$&$0$&$1$&$0$
\\$40_2$&$0$&$1$&$0$&$0$&$0$&$0$&$0$&$0$&$0$&$0$&$1$&$0$&$1$&$0$&$0$
\\$40_2^*$&$0$&$0$&$1$&$0$&$0$&$0$&$0$&$0$&$0$&$0$&$1$&$1$&$0$&$0$&$0$
\\$74$&$1$&$0$&$0$&$0$&$0$&$0$&$0$&$0$&$1$&$1$&$0$&$0$&$0$&$0$&$0$
\\$160$&$0$&$1$&$0$&$0$&$1$&$0$&$0$&$0$&$0$&$1$&$0$&$0$&$0$&$0$&$0$
\\$160^*$&$0$&$0$&$1$&$1$&$0$&$0$&$0$&$0$&$1$&$0$&$0$&$0$&$0$&$0$&$0$
\\$280$&$0$&$0$&$0$&$0$&$0$&$0$&$0$&$1$&$0$&$0$&$0$&$0$&$0$&$0$&$0$
\\$280^*$&$0$&$0$&$0$&$0$&$0$&$0$&$1$&$0$&$0$&$0$&$0$&$0$&$0$&$0$&$0$
\end{tabular}\end{center}

\newpage
\section{\texorpdfstring{$\SU_5(3)$}{SU(5,3)}}

\subsection*{Block information}

\begin{center}\begin{tabular}{ccc}
\hline Block & Kernel & Modules
\\\hline $1$ & $1$ & $1,5,5^*,10,10^*,15,15^*,24,30,30^*,45_1,45_1^*, 45_2,45_2^*,51,65,65^*,105,105^*,$
\\&&$126,126^*,135,135^*,165,165^*,174,174^*,195,195^*,199,315,315^*,375,375^*,390,390^*,$
\\ && $420,420^*,450,450^*,476,510,510^*,720,720^*,835,835^*,855,855^*,1125,1215,1215^*,$
\\&&$1440,1440^*,1545,1545^*,1701,1701^*,2430,2430^*,2535,2535^*,2574,2574^*,3465,3465^*,$
\\&&$3780,3780^*,4199,6810,6810^*,7350,7350^*,8505,8505^*,10824,19260,19260^*,26730,26730^*$
\\ $2$ & $1$ & $59049$
\\\hline\end{tabular}\end{center}

\subsection*{Extension information}

\begin{center}\begin{tabular}{ccc}\hline
Module & Tested to & $\Ext^1=k$ 
\\\hline $1$ & $199$ & $45_1,45_1^*,174,174^*,199$
\\ $5$ & $65$ & $30,45_1,51,65$
\\ $10$ & $65$ & $30$
\\ $15$ & $65$ & $45_1$
\\ $24$ & $65$ & $45_1,45_1^*,51$
\\ $30$ & $65$ & $5,10,30^*$
\\ $45_1$ & $51$ & $5,15,24$
\\ $45_2$ & $51$ & None
\\ $51$ & $51$ & $5,24$
\\ \hline\end{tabular}\end{center}

\newpage
\section{\texorpdfstring{$\SU_5(5)$}{SU(5,5)}}

\subsection*{Block information}

These modules are only necessarily complete up to dimension 1200.

\begin{center}\begin{tabular}{ccc}
\hline Block & Kernel & Modules
\\\hline $1$ & $1$ & $1,5,5^*,10,10^*,15,15^*,23,35,35^*,40,40^*,45,45^*,50,50^*,70_1,70_1^*,70_2,70_2^*,75,$
\\&&$103,103^*,105,105^*,121,121^*,160,160^*,170,170^*,175_1,175_1^*,175_2,175_2^*,175_3,175_3^*,$
\\&&$185,185^*,200,235,235^*,255,255^*,280,280^*,305,305^*,315,315^*,320,320^*,365,365^*,381,$
\\&&$410,410^*,420,420^*,450_1,450_1^*,450_2,450_2^*,535,535^*,540,540^*,545,545^*,615,615^*,670,670^*,$
\\ &&$683,683^*,700_1,700_1^*,700_2,700_2^*,794,875,875^*,895,895^*,955,955^*,1000,1050_1,1050_1^*$
\\ &&$1050_2,1050_2^*,1055,1055^*,1063,1063^*,1090,1090^*1115,1115^*,1150,1150^*,1200,1200^*$,...
\\ $2$ & $1$ & $9765625$
\\\hline\end{tabular}\end{center}

\subsection*{Extension information}

\begin{center}\begin{tabular}{ccc}\hline
Module & Tested to & $\Ext^1=k$ 
\\\hline $1$ & $235$ & $23,185,185^*$
\\ $5$ & $45$ & None
\\ $10$ & $45$ & None
\\ $15$ & $45$ & None
\\ $23$ & $45$ & $1$
\\ $35$ & $45$ & None
\\ $40$ & $45$ & None
\\ $45$ & $45$ & None
\\ \hline\end{tabular}\end{center}

\newpage
\section{\texorpdfstring{$3\cdot\Sp_4(2)'$}{3.Sp(4,2)'}}


\subsection*{Block information}

\begin{center}\begin{tabular}{ccc}
\hline Block & Kernel & Modules
\\ \hline 1 & 3 & $1$, $4_1$, $4_2$
\\ 2 & 3 & $8_1$
\\ 3 & 3 & $8_2$
\\ 4 & 1 & $3_1$, $3_2^*$, $9$
\\ 5 & 1 & Dual to block 4
\\\hline\end{tabular}\end{center}

\subsection*{Identification of modules} 

\begin{enumerate}
\item Fix $3_1$. Then $3_2$ is the image under the field automorphism.
\item $4_1$ is a submodule of $3_1\otimes 3_2$ (and $4_2$ is a quotient).
\item $8_1$ is a submodule of $3_1\otimes 3_1^*$.
\end{enumerate}

\subsection*{Highest weight labels}

\begin{center}\begin{tabular}{ccccc}
\hline $\lambda$ &00&10&01&11
\\\hline $L(\lambda)$&$1$&$4_1$&$4_2$&$8_1\oplus 8_2$
\\ $\dim(W(\lambda))$&$1$&$4$&$4$&$16$
\\ \hline
\end{tabular}\end{center}

\subsection*{Extension information}

\begin{center}\begin{tabular}{c|ccc}
Block 1&$1$&$4_1$&$4_2$
\\ \hline $1$&$0$&$1$&$1$
\\ $4_1$&$1$&$0$&$0$
\\ $4_2$&$1$&$0$&$0$
\end{tabular}\end{center}

\begin{center}\begin{tabular}{c|ccc}
Block 4&$3_1$&$3_2^*$&$9$
\\ \hline $3_1$&$0$&$1$&$1$
\\ $3_2^*$&$1$&$0$&$1$
\\ $9$&$1$&$1$&$0$
\end{tabular}\end{center}

\newpage
\section{\texorpdfstring{$\Sp_4(3)$}{Sp(4,3)}}

\subsection*{Block information}

\begin{center}\begin{tabular}{ccc}
\hline Block & Kernel & Modules
\\\hline1 & 2 & $1$, $5$, $10$, $14$, $25$
\\2 & 2 & $81$
\\3 & 1 & $4$, $16$, $40$
\\\hline\end{tabular}\end{center}

\subsection*{Highest weight labels}

\begin{center}\begin{tabular}{ccccccc}
\hline $\lambda$&00&01&20&02&21&22
\\\hline $L(\lambda)$ &$1$&$5$&$10$&$14$&$25$&$81$
\\ $\dim(W(\lambda))$&$1$&$5$&$10$&$14$&$35$&$81$
\\\hline
\end{tabular}\end{center}

\begin{center}\begin{tabular}{cccc}
\hline $\lambda$&10&11&12
\\\hline $L(\lambda)$ &$4$&$16$&$40$
\\ $\dim(W(\lambda))$&$4$&$16$&$40$
\\\hline
\end{tabular}\end{center}

\subsection*{Extension information}

\begin{center}\begin{tabular}{c|ccccc}
Block 1&$1$&$5$&$10$&$14$&$25$
\\ \hline $1$&$0$&$1$&$0$&$1$&$0$
\\ $5$&$1$&$1$&$1$&$1$&$1$
\\ $10$&$0$&$1$&$1$&$1$&$1$
\\ $14$&$1$&$1$&$1$&$0$&$0$
\\ $25$&$0$&$1$&$1$&$0$&$0$
\end{tabular}\end{center}

\begin{center}\begin{tabular}{c|ccc}
Block 3&$4$&$16$&$40$
\\ \hline $4$&$1$&$2$&$1$
\\ $16$&$2$&$1$&$1$
\\ $40$&$1$&$1$&$0$
\end{tabular}\end{center}

\subsection*{Socle layers of projective modules}

\[
\begin{array}{c}
1
\\5\oplus14
\\1\oplus1\oplus5\oplus10\oplus14
\\1\oplus5\oplus5\oplus10\oplus14\oplus25
\\1\oplus1\oplus5\oplus10\oplus14
\\5\oplus14
\\1
\end{array}\qquad
\begin{array}{c}
5
\\1\oplus5\oplus10\oplus14\oplus25
\\1\oplus5\oplus5\oplus5\oplus5\oplus10\oplus10\oplus14\oplus25
\\1\oplus1\oplus5\oplus5\oplus5\oplus10\oplus10\oplus10\oplus14\oplus14\oplus25\
\oplus25
\\1\oplus5\oplus5\oplus5\oplus5\oplus10\oplus10\oplus14\oplus25
\\1\oplus5\oplus10\oplus14\oplus25
\\5
\end{array}
\]
\[ \begin{array}{c}
10
\\10\oplus25
\\10\oplus25
\\10\oplus10
\\10\oplus25
\\10\oplus25
\\10\oplus10
\\5\oplus10\oplus14\oplus25
\\1\oplus5\oplus5\oplus10\oplus10\oplus10\oplus14\oplus25
\\1\oplus5\oplus5\oplus5\oplus10\oplus10\oplus10\oplus14\oplus14\oplus25
\\1\oplus5\oplus5\oplus10\oplus10\oplus10\oplus14\oplus25
\\5\oplus10\oplus14\oplus25
\\10
\end{array}\qquad
\begin{array}{c}
14
\\1\oplus5\oplus10
\\1\oplus5\oplus10\oplus14\oplus14\oplus25
\\1\oplus5\oplus5\oplus10\oplus10\oplus14
\\1\oplus5\oplus10\oplus14\oplus14\oplus25
\\1\oplus5\oplus10
\\14
\end{array}\qquad
\begin{array}{c}
25
\\10
\\10
\\25
\\10
\\10
\\25
\\5\oplus10
\\5\oplus10\oplus14\oplus25
\\1\oplus5\oplus5\oplus10\oplus25
\\5\oplus10\oplus14\oplus25
\\5\oplus10
\\25
\end{array}\]

Remark: these socle layers also appear in David Benson's paper \emph{Projective modules for the group of twenty-seven lines on a cubic surface} (p.1019). There are minor errors for $10$ and $25$.

\newpage
\section{\texorpdfstring{$\Sp_4(5)$}{Sp(4,5)}}

\subsection*{Block information}

\begin{center}\begin{tabular}{ccc}
\hline Block & Kernel & Modules
\\\hline 1&2&$1$, $5$, $10$, $13$, $30$, $35_1$, $35_2$, $55$, $68$, $86$, $105$, $115$, $260$, $355$
\\ 2&2&$625$
\\ 3&1&$4$, $12$, $20$, $40$, $52$, $80$, $140_1$, $140_2$, $204$, $420$
\\\hline\end{tabular}\end{center}

\subsection*{Highest weight labels}

\begin{center}\begin{tabular}{cccccccccccccccc}
\hline $\lambda$&00&01&20&02&03&21&40&04&22&23&41&42&24&43&44
\\\hline $L(\lambda)$ &$1$&$5$&$10$&$13$&$30$&$35_1$&$35_2$&$55$&$68$&$86$&$105$&$115$&$260$&$355$&$625$
\\ $\dim(W(\lambda))$&$1$&$5$&$10$&$14$&$30$&$35$&$35$&$55$&$81$&$154$&$105$&$220$&$260$&$390$&$625$
\\\hline
\end{tabular}\end{center}

\begin{center}\begin{tabular}{ccccccccccc}
\hline $\lambda$&10&11&30&12&31&13&32&14&33&34
\\\hline $L(\lambda)$ &$4$&$12$&$20$&$40$&$52$&$80$&$140_1$&$140_2$&$204$&$420$
\\ $\dim(W(\lambda))$&$4$&$16$&$20$&$40$&$64$&$80$&$140$&$140$&$256$&$420$
\\\hline
\end{tabular}\end{center}

\subsection*{Extension information}

\begin{center}\begin{tabular}{c|cccccccccccccc}
Block 1&$1$&$5$&$10$&$13$&$30$&$35_1$&$35_2$&$55$&$68$&$86$&$105$&$115$&$260$&$355$
\\\hline  $1$&$0$&$0$&$0$&$1$&$0$&$0$&$0$&$0$&$1$&$0$&$0$&$1$&$0$&$0$
\\ $5$&$0$&$0$&$0$&$0$&$0$&$0$&$0$&$0$&$1$&$1$&$0$&$0$&$0$&$1$
\\ $10$&$0$&$0$&$0$&$0$&$1$&$0$&$0$&$1$&$1$&$0$&$0$&$1$&$0$&$1$
\\ $13$&$1$&$0$&$0$&$1$&$0$&$1$&$0$&$0$&$1$&$1$&$0$&$0$&$1$&$0$
\\ $30$&$0$&$0$&$1$&$0$&$0$&$1$&$1$&$0$&$0$&$0$&$0$&$0$&$1$&$0$
\\ $35_1$&$0$&$0$&$0$&$1$&$1$&$1$&$0$&$1$&$0$&$1$&$0$&$0$&$1$&$1$
\\ $35_2$&$0$&$0$&$0$&$0$&$1$&$0$&$0$&$1$&$1$&$0$&$0$&$0$&$0$&$1$
\\ $55$&$0$&$0$&$1$&$0$&$0$&$1$&$1$&$0$&$0$&$0$&$0$&$0$&$0$&$0$
\\ $68$&$1$&$1$&$1$&$1$&$0$&$0$&$1$&$0$&$1$&$1$&$1$&$0$&$0$&$0$
\\$86$&$0$&$1$&$0$&$1$&$0$&$1$&$0$&$0$&$1$&$0$&$0$&$0$&$0$&$0$
\\$105$&$0$&$0$&$0$&$0$&$0$&$0$&$0$&$0$&$1$&$0$&$1$&$1$&$0$&$0$
\\ $115$&$1$&$0$&$1$&$0$&$0$&$0$&$0$&$0$&$0$&$0$&$1$&$0$&$0$&$0$
\\ $260$&$0$&$0$&$0$&$1$&$1$&$1$&$0$&$0$&$0$&$0$&$0$&$0$&$0$&$0$
\\ $355$&$0$&$1$&$1$&$0$&$0$&$1$&$1$&$0$&$0$&$0$&$0$&$0$&$0$&$0$
\end{tabular}\end{center}

\begin{center}\begin{tabular}{c|cccccccccc}
Block 3 &$4$&$12$&$20$&$40$&$52$&$80$&$140_1$&$140_2$&$204$&$420$
\\\hline $4$&$0$&$1$&$0$&$0$&$1$&$1$&$1$&$0$&$1$&$0$
\\$12$&$1$&$1$&$0$&$1$&$1$&$0$&$0$&$1$&$1$&$0$
\\$20$&$0$&$0$&$0$&$1$&$0$&$2$&$0$&$1$&$0$&$0$
\\$40$&$0$&$1$&$1$&$1$&$0$&$0$&$0$&$1$&$1$&$1$
\\$52$&$1$&$1$&$0$&$0$&$1$&$1$&$1$&$0$&$1$&$0$
\\$80$&$1$&$0$&$2$&$0$&$1$&$0$&$0$&$0$&$0$&$1$
\\$140_1$&$1$&$0$&$0$&$0$&$1$&$0$&$1$&$0$&$0$&$0$
\\$140_2$&$0$&$1$&$1$&$1$&$0$&$0$&$0$&$0$&$0$&$0$
\\$204$&$1$&$1$&$0$&$1$&$1$&$0$&$0$&$0$&$0$&$0$
\\$420$&$0$&$0$&$0$&$1$&$0$&$1$&$0$&$0$&$0$&$0$
\end{tabular}\end{center}

\newpage
\section{\texorpdfstring{$\Sp_4(7)$}{Sp(4,7)}}

\subsection*{Block information}

\begin{center}\begin{tabular}{ccc}
\hline Block & Kernel & Modules
\\\hline 1 & 2 & $1$, $5$, $10$, $14$, $25$, $35_1$, $35_2$, $54$, $71$, $84$, $91$, $105$, $140_1$, $149$, $154$, $199$,
\\ && $206$, $231$, $260$, $315$, $365$, $455$, $595$, $786$, $959$, $1330$, $1645$
\\ 2 & 2 & $2401$
\\ 3 & 1 & $4$, $16$, $20$, $24$, $44$, $56$, $76$, $116$, $140_2$, $140_3$, $180$,
\\&&$224$, $300$, $336$, $420$, $460$, $560$, $580$, $924$, $1180$, $1820$
\\\hline\end{tabular}\end{center}

\subsection*{Highest weight labels}

\begin{center}\begin{tabular}{ccccccccccccccc}
\hline $\lambda$&00&01&20&02&03&21&40&04&22&60&05&41&06&42
\\\hline $L(\lambda)$ &$1$&$5$&$10$&$14$&$25$&$35_1$&$35_2$&$54$&$71$&$84$&$91$&$105$&$140_1$&$149$
\\ $\dim(W(\lambda))$&$1$&$5$&$10$&$14$&$30$&$35$&$35$&$55$&$81$&$84$&$91$&$105$&$140$&$220$
\\\hline $\lambda$&23&25&24&61&44&63&43&62&26&45&64&46&65&66
\\\hline $L(\lambda)$ &$154$&$199$& $206$&$231$&$260$&$315$&$365$&$455$&$595$&$786$&$959$&$1330$&$1645$&$2401$
\\ $\dim(W(\lambda))$&$154$&$405$&$260$&$231$&$625$&$770$&$390$&$455$&$595$&$935$&$1190$&$1330$&$1729$&$2401$
\\\hline
\end{tabular}\end{center}

\begin{center}\begin{tabular}{cccccccccccc}
\hline $\lambda$&10&11&30&12&31&50&13&51&32&14&33
\\\hline $L(\lambda)$ & $4$&$16$&$20$&$24$&$44$&$56$&$76$&$116$&$140_2$&$140_3$&$180$
\\ $\dim(W(\lambda))$&$4$&$16$&$20$&$40$&$64$&$56$&$80$&$160$&$140$&$140$&$256$
\\ \hline
 $\lambda$&15&52&16&34&35&53&54&36&55&56
\\ \hline $L(\lambda)$ & $224$&$300$&$336$&$420$&$460$&$560$&$580$&$924$&$1180$&$1820$
\\ $\dim(W(\lambda))$&$224$&$324$&$336$&$420$&$640$&$560$&$880$&$924$&$1296$&$1820$
\\ \hline
\end{tabular}\end{center}

\subsection*{Extension information}

\begin{center}
\begin{tabular}{c|c}
Block 1 &  $\Ext^1=k$
\\ \hline $1$ & $54$, $149$, $260$
\\ $5$&$25$, $149$, $199$, $315$, $365$, $786$
\\ $10$&$71$, $149$, $206$, $260$, $365$, $786$, $959$
\\ $14$&$260$, $365$, $959$
\\ $25$&$5$, $25$, $71$, $260$, $365$, $786$
\\ $35_1$ & $206$, $315$, $365$, $786$, $959$, $1645$
\\ $35_2$ & $91$, $140_1$, $206$, $959$, $1645$
\\ $54$ & $1$, $71$, $105$, $206$, $786$, $1330$
\\ $71$ & $10$, $25$, $54$, $71$, $149$, $154$, $199$, $260$, $595$, $786$
\\ $84$ & $91$, $140_1$, $206$, $1645$
\\ $91$ & $35_2$, $84$, $105$, $1330$
\\ $105$& $54$, $91$, $140_1$, $154$, $199$, $595$, $1645$
\\ $140_1$ & $35_2$, $84$, $105$
\\ $149$ & $1$, $5$, $10$, $71$, $149$, $206$, $365$, $786$
\\ $154$ & $71$, $105$, $154$, $595$, $786$, $1330$
\\ $199$ & $5$, $71$, $105$, $206$
\\ $206$ & $10$, $35_1$, $35_2$, $54$, $84$, $149$, $199$, $231$
\\ $231$ & $206$, $365$, $959$
\\ $260$ & $1$, $10$, $14$, $25$, $71$, $365$
\\ $315$ & $5$, $35_1$, $455$
\\ $365$ & $5$, $10$, $14$, $25$, $35_1$, $149$, $231$, $260$, $365$, $455$
\\ $455$ & $315$, $365$, $455$
\\ $595$ & $71$, $105$, $154$
\\ $786$ & $5$, $10$, $25$, $35_1$, $54$, $71$, $149$, $154$
\\ $959$ & $10$, $14$, $35_1$, $35_2$, $231$
\\ $1330$ & $54$, $91$, $154$
\\ $1645$ & $35_1$, $35_2$, $84$, $105$
\end{tabular}
\end{center}

\begin{center}
\begin{tabular}{c|cc}
Block 3 & $\Ext^1=k$ & $\Ext^1=k^2$
\\ \hline $4$ & $76$, $116$, $180$, $300$, $460$, $580$ &
\\ $16$ & $24$, $180$, $300$, $420$, $460$, $560$, $580$, $1180$ &
\\ $20$ & $44$, $180$, $224$, $420$, $580$, $1180$ &
\\ $24$ & $16$, $24$, $76$, $300$, $460$, $580$ &
\\ $44$ & $20$, $76$, $116$, $140_3$, $336$, $460$ &
\\ $56$ & $140_3$, $336$ & $224$
\\ $76$ & $4$, $24$, $44$, $76$, $140_2$, $180$, $460$, $580$, $924$, $1180$ &
\\ $116$ & $4$, $44$, $180$, $224$, $300$, $420$, $1180$ &
\\ $140_2$ & $76$, $140_2$, $140_3$, $336$, $460$, $924$ &
\\ $140_3$ & $44$, $56$, $140_2$, $924$, $1180$, $1820$ &
\\ $180$ & $4$, $16$, $20$, $76$, $116$, $180$, $300$, $460$ &
\\ $224$ & $20$, $116$, $1820$ & $56$
\\ $300$ & $4$, $16$, $24$, $116$, $180$, $300$, $420$, $560$, $580$ &
\\ $336$ & $44$, $56$, $140_2$ &
\\ $420$ & $16$, $20$, $116$, $300$ &
\\ $460$ & $4$, $16$, $24$, $44$, $76$, $140_2$, $180$ &
\\ $560$ & $16$, $300$, $560$ &
\\ $580$ & $4$, $16$, $20$, $24$, $76$, $300$ &
\\ $924$ & $76$, $140_2$, $140_3$ &
\\ $1180$ & $16$, $20$, $76$, $116$, $140_3$ &
\\ $1820$ & $140_3$, $224$ &
\end{tabular}
\end{center}

\newpage
\section{\texorpdfstring{$\Sp_6(2)$}{Sp(6,2)}}

\subsection*{Block information}

\begin{center}\begin{tabular}{ccc}
\hline Block & Kernel & Modules
\\\hline  1 & $1$ & $1$, $6$, $8$, $14$, $48$, $64$, $112$
\\ 2 & $1$ & $512$
\\\hline\end{tabular}\end{center}

\subsection*{Highest weight labels}

\begin{center}\begin{tabular}{ccccccccc}
\hline $\lambda$&000&100&001&010&101&110&011&111
\\\hline $L(\lambda)$ &$1$&$6$&$8$&$14$&$48$&$64$&$112$&$512$
\\ $\dim(W(\lambda))$&$1$&$6$&$14$&$14$&$70$&$64$&$126$&$512$
\\\hline
\end{tabular}\end{center}

\subsection*{Extension information}

\begin{center}\begin{tabular}{c|ccccccc}
Block 1& $1$&$6$&$8$&$14$&$48$&$64$&$112$
\\ \hline $1$&$0$&$1$&$0$&$0$&$1$&$0$&$0$
\\ $6$&$1$&$0$&$1$&$2$&$0$&$0$&$0$
\\ $8$&$0$&$1$&$0$&$0$&$0$&$0$&$1$
\\ $14$&$0$&$2$&$0$&$0$&$0$&$1$&$0$
\\ $48$&$1$&$0$&$0$&$0$&$0$&$0$&$1$
\\ $64$&$0$&$0$&$0$&$1$&$0$&$1$&$0$
\\ $112$&$0$&$0$&$1$&$0$&$1$&$0$&$0$
\end{tabular}\end{center}

\subsection*{Tensor products with the Steinberg}

\begin{center}\begin{tabular}{c|c}
$M$ & $M\otimes 512$
\\\hline $6$ & $P(112)\oplus 512$
\\ $8$ & $P(64)\oplus 512$
\\ $14$ & $P(48)\oplus 512^{\oplus 3}$
\\ $48$ & $P(14)\oplus P(64)\oplus P(112)\oplus 512^{\oplus 7}$
\\ $64$ & $P(8)\oplus P(48)\oplus P(64)^{\oplus 2}\oplus P(112)^{\oplus 3}\oplus 512^{\oplus 11}$
\\ $112$ & $P(6)\oplus P(48)\oplus P(64)^{\oplus 3}\oplus P(112)^{\oplus 4}\oplus 512^{\oplus 19}$
\end{tabular}\end{center}

\subsection*{Symmetric and exterior squares of modules}

\begin{center}\begin{tabular}{c|ccc}
Block 1&$6$&$8$&$14$
\\\hline $S^2(M)$ &$1,14/6$&$1/6/14/6/1,8$&$1\oplus (14/6/1/48/1/6/14,14)$
\\ $\Lambda^2(M)$&$1\oplus 14$&$1/6/14/6/1$&$1\oplus (14/6/1/48/1/6/14)$
\end{tabular}\end{center}

\subsection*{Actions of unipotent elements on modules}

\begin{center}\begin{tabular}{c|ccccccc}
Block 1 & $1$&$6$&$8$&$14$&$48$&$64$&$112$
\\\hline 2A&$1$&$2,1^4$&$2^4$&$2^4,1^6$&$2^{24}$&$2^{24},1^{16}$&$2^{56}$
\\ 2B&$1$&$2^2,1^2$&$2^2,1^4$&$2^6,1^2$&$2^{20},1^8$&$2^{32}$&$2^{52},1^8$
\\ 2C&$1$&$2^2,1^2$&$2^4$&$2^6,1^2$&$2^{24}$&$2^{32}$&$2^{56}$
\\ 2D&$1$&$2^3$&$2^4$&$2^6,1^2$&$2^{24}$&$2^{32}$&$2^{56}$
\\ 4A[B],4D[B]&$1$&$3^2$&$3^2,1^2$&$4^2,3^2$&$4^8,3^4,1^4$&$4^{16}$&$4^{24},3^4,1^4$
\\ 4B[C],4C[C]&$1$&$4,1^2$&$4^2$&$4^3,2$&$4^{12}$&$4^{16}$&$4^{28}$
\\ 4E[C]&$1$&$4,2$&$4^2$&$4^3,2$&$4^{12}$&$4^{16}$&$4^{28}$
\\ 8A[A], 8B[D]&$1$&$6$&$6,2$&$8,6$&$8^4,6^2,2^2$&$8^8$&$8^{12},6^2,2^2$
\end{tabular}\end{center}

\subsection*{Traces of semisimple elements on modules}

\begin{center}\begin{tabular}{c|cccccccc}
Block 1,2 & $1$&$6$&$8$&$14$&$48$&$64$&$112$&$512$
\\\hline 3A&$1$&$3$&$-4$&$2$&$-12$&$4$&$-8$&$-16$
\\3B&$1$&$-3$&$-1$&$5$&$3$&$-8$&$-5$&$8$
\\3C&$1$&$0$&$2$&$-1$&$0$&$-2$&$-2$&$-4$
\\5A&$1$&$1$&$-2$&$-1$&$-2$&$-1$&$2$&$2$
\\7A&$1$&$-1$&$1$&$0$&$-1$&$1$&$0$&$1$
\\9A[B]&$1$&$0$&$-1$&$-1$&$0$&$1$&$1$&$-1$
\\15A[A,A]&$1$&$-2$&$1$&$2$&$-2$&$-1$&$2$&$-1$
\end{tabular}\end{center}

\newpage
\section{\texorpdfstring{$\Sp_6(3)$}{Sp(6,3)}}

\subsection*{Block information}

\begin{center}\begin{tabular}{ccc}
\hline Block & Kernel & Modules
\\ \hline Block 1 & $2$ & $1$, $13$, $21$, $57$, $84$, $90$, $189$, $358$, $594$, $623$, $903$, $1302$, $2457$, $4655$
\\ Block 2 & $2$ & $19683$
\\ Block 3 & $1$ & $6$, $14$, $50$, $126$, $216$, $286$, $316$, $378$, $1260$, $1610$, $3450$, $8316$
\\\hline\end{tabular}\end{center}

\subsection*{Highest weight labels}

\begin{center}\begin{tabular}{ccccccccc}
\hline $\lambda$&000&010&200&101&002&020&210&111
\\\hline $L(\lambda)$&$1$&$13$&$21$&$57$&$84$&$90$&$189$&$358$
\\ $\dim(W(\lambda))$&$1$&$14$&$21$&$70$&$84$&$90$&$189$&$512$
\\\hline
$\lambda$&012&202&220&121&022&212&222&
\\\hline $L(\lambda)$&$594$&$623$&$903$&$1302$&$2457$&$4655$&$19683$&
\\ $\dim(W(\lambda))$&$594$&$1078$&$924$&$2205$&$2457$&$5720$&$19683$&
\\ \hline
\end{tabular}\end{center}

\begin{center}\begin{tabular}{ccccccc}
\hline $\lambda$&100&001&110&011&201&120
\\\hline $L(\lambda)$&$6$&$14$&$50$&$126$&$216$&$286$
\\ $\dim(W(\lambda))$&$6$&$14$&$64$&$126$&$216$&$350$
\\\hline
$\lambda$&021&102&211&112&221&122
\\\hline $L(\lambda)$&$316$&$378$&$1260$&$1610$&$3450$&$8316$
\\ $\dim(W(\lambda))$&$616$&$378$&$1386$&$2240$&$5460$&$8316$
\\ \hline
\end{tabular}\end{center}

\subsection*{Extension information}

\begin{center}\begin{tabular}{c|cccccccccccccc}
Module & $1$&$13$&$21$&$57$&$84$&$90$&$189$&$358$&$594$&$623$&$903$&$1302$&$2457$&$4655$
\\ \hline $1$&$0$&$1$&$0$&$0$&$0$&$1$&$0$&$0$&$0$&$0$&$0$&$0$&$0$&$1$
\\ $13$&$1$&$1$&$0$&$1$&$0$&$0$&$0$&$2$&$0$&$1$
\\ $21$&$0$&$0$&$0$&$1$&$1$&$1$&$0$&$0$
\\ $57$&$0$&$1$&$1$&$0$&$0$&$1$&$1$
\\ $84$&$0$&$0$&$1$&$0$&$0$&$1$
\\ $90$&$1$&$0$&$1$&$1$&$1$&$1$
\\ $189$&$0$&$0$&$0$&$1$
\\ $358$&$0$&$2$&$0$&
\\ $594$&$0$&$0$&
\\ $623$&$0$&$1$& 
\\ $903$&$0$& 
\\ $1302$&$0$&
\\ $2457$&$0$&
\\ $4655$&$1$&
\end{tabular}\end{center}

\begin{center}\begin{tabular}{c|cccccccccccc}
Module & $6$&$14$&$50$&$126$&$216$&$286$&$316$&$378$&$1260$&$1610$&$3450$&$8316$
\\ \hline $6$&$0$&$0$&$1$&$1$&$0$&$0$&$0$&$0$&$1$&$1$&&$0$
\\ $14$&$0$&$1$&$1$&$0$&$0$&$1$&&&&&&$0$
\\ $50$&$1$&$1$&$0$&$2$&&&&&&&&$0$
\\ $126$&$1$&$0$&$2$&$1$&&&&&&&&$0$
\\ $216$&$0$&$0$&&&&&&&&&&$0$
\\ $286$&$0$&$1$&&&&&&&&&&$0$
\\ $316$&$0$&&&&&&&&&&&$1$
\\ $378$&$0$&&&&&&&&&&&$1$
\\ $1260$&$1$&&&&&&&&&&&$0$
\\ $1610$&$1$&&&&&&&&&&&$0$
\\ $3450$&&&&&&&&&&&&$0$
\\ $8316$&$0$&$0$&$0$&$0$&$0$&$0$&$1$&$1$&$0$&$0$&$0$&$0$
\end{tabular}\end{center}

\subsection*{Actions of unipotent elements on modules}

\begin{small}
\begin{center}\begin{tabular}{c|cccccc}
Block 1&$1$&$13$&$21$&$57$&$84$&$90$
\\\hline 3A[2]&$1$&$2^4,1^5$&$3,2^4,1^{10}$&$3^5,2^{16},1^{10}$&$3^{14},2^{16},1^{10}$&$3^{10},2^{20},1^{20}$
\\ 3C,3D&$1$&$3,2^4,1^2$&$3^3,2^4,1^4$&$3^{14},2^4,1^7$&$3^{23},2^4,1^7$&$3^{22},2^8,1^8$
\\ 3E[2]&$1$&$3^3,1^4$&$3^6,1^3$&$3^{18},1^3$&$3^{27},1^3$&$3^{27},1^9$
\\ 3G&$1$&$3^4,1$&$3^7$&$3^{19}$&$3^{27},1^3$&$3^{30}$
\\ 9A[A][2]&$1$&$5,4^2$&$7,4^2,3,1^3$&$8^2,7,6^2,5^3,3,2^2$&$9^3,8^2,7,6^2,5^3,3,2^2$&$9,8^2,7^3,6^2,5^2,4^2,3^3,2^2,1$
\\9C[A][2],9E[A][2]&$1$&$5^2,3$&$7,5,3^3$&$9,7^2,6^4,3^3,1$&$9^4,7^2,6^4,3^3,1$&$9^3,7,6^4,5^3,3^5,1^2$
\\9G[E][2]&$1$&$9,3,1$&$9^2,3$&$9^6,3$&$9^9,3$&$9^9,3^3$
\end{tabular}\end{center}

\begin{center}\begin{tabular}{c|ccccc}
Block 1&$189$&$358$&$594$&$623$&$903$
\\\hline 3A[2]&$3^{24},2^{36},1^{45}$&$3^{77},2^{44},1^{39}$&$3^{159},2^{36},1^{45}$&$3^{189},2^{20},1^{16}$&$3^{199},2^{100},1^{106}$
\\ 3C,3D&$3^{51},2^{12},1^{12}$&$3^{113},2^8,1^3$&$3^{186},2^{12},1^{12}$&$3^{199},2^8,1^{10}$&$3^{297},2^4,1^4$
\\ 3E[2]&$3^{63}$&$3^{117},1^7$&$3^{198}$&$3^{201},2^8,1^4$&$3^{297},2^6$
\\ 3G&$3^{63}$&$3^{118},1^4$&$3^{198}$&$3^{206},1^5$&$3^{301}$
\\ 9A[A][2]&$9^77^36^65^34^63^91^3$&$9^{24}8^276^{12}54^23^{10}2^2$&$9^{52}7^36^65^34^63^91^3$&$9^{63}6^25^54^43$&$9^{65}7^46^{30}4^23^{34}$
\\9C[A][2],9E[A][2]&$9^8,6^{12},3^{15}$&$9^{25},7^2,6^{12},5^2,3^{12},1$&$9^{53},6^{12},3^{15}$&$9^{63},6^4,5^4,3^4$&$9^{66},7,6^{32},5,3^{35}$
\\9G[E][2]&$9^{21}$&$9^{39},3^2,1$&$9^{66}$&$9^{67},6^2,5,3$&$9^{99},6^2$
\end{tabular}\end{center}

\begin{center}\begin{tabular}{c|ccccc}
Block 3&$6$&$14$&$50$&$126$&$216$
\\\hline 3A[2]&$2,1^4$&$2^5,1^4$&$3^4,2^{11},1^{16}$&$3^{16},2^{29},1^{20}$&$3^{30},2^{45},1^{36}$
\\3C,3D&$2^2,1^2$&$3^2,2^2,1^4$&$3^{10},2^8,1^4$&$3^{32},2^{10},1^{10}$&$3^{60},2^{12},1^{12}$
\\ 3E[2]&$2^3$&$3^2,2^4$&$3^{14},2^4$&$3^{36},2^9$&$3^{72}$
\\ 3G&$3^2$&$3^4,1^2$&$3^{16},1^2$&$3^{42}$&$3^{72}$
\\9A[A][2]&$4,1^2$&$5^2,4$&$8,7^2,6,4^3,3^2,2,1^2$&$9^2,8^4,7^2,6^4,5^4,4,3^2,2^4$&$9^8,8^3,6^9,5^6,4^3,3^6,2^3$
\\9C[A][2],9E[A][2]&$4,2$&$6,4^2$&$8^2,6^3,3^4,2^2$&$9^4,8^2,6^8,4^3,3^4,2$&$9^{10},6^{15},3^{12}$
\\ 9G[E][2]&$6$&$8,6$&$9^4,8,6$&$9^{12},6^3$&$9^{24}$
\end{tabular}\end{center}

\begin{center}\begin{tabular}{c|ccc}
Block 3&$286$&$316$&$378$
\\\hline 3A[2]&$3^{56},2^{39},1^{40}$&$3^{90},2^{15},1^{16}$&$3^{84},2^{45},1^{36}$
\\ 3C,3D&$3^{92},2^2,1^6$&$3^{98},2^8,1^6$&$3^{114},2^{12},1^{12}$
\\ 3E[2]&$3^{92},2^4,1^2$&$3^{102},2^4,1^2$&$3^{126}$
\\ 3G&$3^{94},1^4$&$3^{104},1^4$&$3^{126}$
\\9A[A][2]&$9^{18},7^2,6^{10},5^2,4,3^{12}$&$9^{28},8,7^4,6,4^3,3^2,2,1^2$&$9^{26},8^3,6^9,5^6,4^3,3^6,2^3$
\\9C[A][2],9E[A][2]&$9^{18},8,6^{12},4^2,3^{12}$&$9^{28},8^3,6^4,3^4,2^2$&$9^{28},6^{15},3^{12}$
\\ 9G[E][2]&$9^{30},8,6,2$&$9^{34},6,4$&$9^{42}$
\end{tabular}\end{center}
\end{small}

\newpage
\section{\texorpdfstring{$\Sp_6(5)$}{Sp(6,5)}}

\subsection*{Block information}

These modules are only necessarily complete up to dimension 2000.

\begin{center}\begin{tabular}{ccc}
\hline Block & Kernel & Modules
\\ \hline Block 1 & 2 & $1$, $14_1$, $21$, $63$, $70$, $90$, $126$, $189$, $295$, $309$, $426$, $525$, $573$, $798$, $813$, $938$, $951_1$, $951_2$,
\\ && $1197$, $1506$, $2135$, $2619$, $3513$, $3585$, $3900$, $4143$, $7314$, $9450$, $12375$, \dots
\\ Block 2 & 2 & $1953125$
\\ Block 3 & 1 & $6$, $14_2$, $56$, $62$, $64$, $202$, $246$, $316$, $330$, $350$, $552$, $666$, $924$, $1100$, $1246$, $1288$, $1644$,
\\ && $1686$, $2240$, $3192$, $3220$, $3794$, $3962$, $4858$, $6300$, $6390$, $6874$, $7686$, \dots
\\\hline\end{tabular}\end{center}


\subsection*{Highest weight labels}

\begin{center}\begin{tabular}{cccccccccccc}
\hline $\lambda$&000&010&200&002&101&020&400&210&030&111&220
\\\hline $L(\lambda)$&$1$&$14_1$&$21$&$63$&$70$&$90$&$126$&$189$&$295$&$309$&$426$
\\ $\dim(W(\lambda))$&$1$&$14$&$21$&$84$&$70$&$90$&$126$&$189$&$385$&$512$&$924$
\\\hline
$\lambda$&301&012&410&103&004&202&040&311&022&121&
\\\hline $L(\lambda)$&$525$&$573$&$798$&$813$&$938$&$951_1$&$951_2$&$1197$&$1506$&$2135$&
\\ $\dim(W(\lambda))$&$525$&$594$&$924$&$1386$&$1001$&$1078$&$1274$&$3072$&$2457$&$2205$&
\\ \hline
\end{tabular}\end{center}

\begin{center}\begin{tabular}{cccccccccc}
\hline $\lambda$&100&001&300&011&110&201&310&102&003
\\\hline $L(\lambda)$&$6$&$14_2$&$56$&$62$&$64$&$202$&$246$&$316$&$330$
\\ $\dim(W(\lambda))$&$6$&$14$&$56$&$126$&$64$&$216$&$448$&$378$&$330$
\\\hline
$\lambda$&120&021&130&211&401&031&302&320&013
\\\hline $L(\lambda)$&$350$&$552$&$666$&$924$&$1100$&$1246$&$1288$&$1644$&$1686$
\\ $\dim(W(\lambda))$&$350$&$616$&$1344$&$1386$&$1100$&$2184$&$2464$&$2016$&$2002$
\\ \hline
\end{tabular}\end{center}

\subsection*{Extension information}

\begin{center}\begin{tabular}{ccc}\hline
Module & Tested to & $\Ext^1=k$ 
\\\hline $1$ & $2619$ & $426$, $951_1$
\\ $14_1$ & $426$ & $309$, $426$
\\ $21$ & $189$ & $63$
\\ $63$ & $70$ & $21$, $63$
\\ \hline\end{tabular}\end{center}

\newpage
\section{\texorpdfstring{$\Sp_6(7)$}{Sp(6,7)}}

\subsection*{Block information}

These modules are only necessarily complete up to dimension 2000.
\begin{center}\begin{tabular}{ccc}
\hline Block & Kernel & Modules
\\ \hline 1 & 2 & 
$1$, $14_1$, $21$, $70$, $84$, $89$, $126_1$, $171$, $189$, $385$, $423$, $462$, $525$, $861$, $875$, $924_1$, $924_2$, $1078$,
\\ && $1251$, $1274$, $1533$, $2079$, $2205$, $2254$, $2388$, $2649$, $2667$, $3003$, $3276$, $3507$, $4830$,
\\ && $5334$, $5922$, $6041$, $6055$, $6087$, $6776$, $8547$, $11914$, $11319$, $12882$, $15092$, and $132$ more
\\ 2 & 2 & $40353607$
\\ $3$ & 1 &
$6$, $14_2$, $56$, $58$, $126_2$, $158$, $172$, $252$, $350$, $378$, $448$, $616$, $786$, $840$, $854$, $942$,
$1344$,
\\ && $1386$, $1544$, $1876$, $2016$, $2086$, $2414$, $3640$, $3696$, $4116$, $4844_1$, $4844_2$, $5950$, $5964$, $6006$,
\\ && $8918$, and $136$ more
\\\hline\end{tabular}\end{center}

\subsection*{Highest weight labels}

\begin{center}\begin{tabular}{cccccccccccc}
\hline $\lambda$&000&010&200&101&002&020&400&012&210&030&111
\\ \hline $L(\lambda)$&$1$&$14_1$&$21$&$70$&$84$&$89$&$126_1$&$171$&$189$&$385$&$423$
\\ $\dim(W(\lambda))$&$1$&$14$&$21$&$70$&$84$&$90$&$126$&$594$&$189$&$385$&$512$
\\ \hline $\lambda$&600&301&103&004&410&220&202&420&040&022&501
\\ \hline $L(\lambda)$&$462$&$525$&$861$&$875$&$924_1$&$924_2$&$1078$&$1251$&$1274$&$1533$&$2079$
\\ $\dim(W(\lambda))$&$462$&$525$&$1386$&$1001$&$924$&$924$&$1078$&$3900$&$1274$&$2457$&$2079$
\\ \hline
\end{tabular}\end{center}

\begin{center}\begin{tabular}{ccccccccccccccc}
\hline $\lambda$&100&001&300&110&011&201&500&120&310&021&510&031&112&401
\\ \hline $L(\lambda)$&$6$&$14_2$&$56$&$58$&$126_2$&$158$&$252$&$350$&$448$&$616$&$786$&$840$&$854$&$942$
\\ $\dim(W(\lambda))$&$6$&$14$&$56$&$64$&$126$&$216$&$252$&$350$&$448$&$616$&$1728$&$2184$&$2240$&$1100$
\\ \hline
\end{tabular}\end{center}

\subsection*{Extension information}

The trivial module only has $1$-cohomology with $89$ and $1251$ for all modules up to $6776$ listed above.

\begin{center}\begin{tabular}{c|ccccccccc}
Block 1&$1$&$14_1$&$21$&$70$&$84$&$89$&$126_1$&$171$&$189$ 
\\\hline $1$&$0$&$0$&$0$&$0$&$0$&$1$&$0$&$0$&$0$
\\$14_1$&$0$&$0$&$0$&$0$&$0$&$0$&$0$&$0$&$0$
\\$21$&$0$&$0$&$0$&$0$&$0$&$0$&$0$&$0$&$0$
\\$70$&$0$&$0$&$0$&$0$&$0$&$0$&$0$&$0$&$0$
\\$84$&$0$&$0$&$0$&$0$&$0$&$0$&$0$&$0$
\\$89$&$1$&$0$&$0$&$0$&$0$&$0$&$0$
\\$126_1$&$0$&$0$&$0$&$0$&$0$&$0$&$0$
\end{tabular}\end{center}

\begin{center}\begin{tabular}{c|ccccc}
Block 3&$6$&$14_2$&$56$&$58$&$126_2$
\\\hline $6$&$0$&$0$&$0$&$1$&$0$
\\$14_2$&$0$&$0$&$0$&$0$&$0$	
\\$56$&$0$&$0$&$0$&$0$&$0$
\\$58$&$1$&$0$&$0$&$0$&$0$
\\$126$&$0$&$0$&$0$&$0$&$0$
\end{tabular}\end{center}

\newpage
\section{\texorpdfstring{$\Sp_8(2)$}{Sp(8,2)}}

\subsection*{Block information}

\begin{center}\begin{tabular}{ccc}
\hline Block & Kernel & Modules
\\\hline $1$ & $1$ & $1,8,16,26,48,128,160,246,416,768,784,2560,3936,4096,12544$
\\ $2$ & $1$ & $65536$
\\\hline\end{tabular}\end{center}

\subsection*{Highest weight labels}

\begin{center}\begin{tabular}{ccccccccc}
\hline $\lambda$&0000&1000&0001&0100&0010&1001&1100&1010
\\\hline $L(\lambda)$ & $1$&$8$&$16$&$26$&$48$&$128$&$160$&$246$
\\ $\dim(W(\lambda))$&$1$&8&42&27&48&288&160&315
\\\hline $\lambda$&0101&0011&0110&1101&1011&1110&0111&1111
\\\hline $L(\lambda)$&$416$&$768$&$784$&$2560$&$3936$&$4096$&$12544$&$65536$
\\ $\dim(W(\lambda))$&792&1056&792&4200&6237&4096&13728&65536
\\ \hline
\end{tabular}\end{center}

\subsection*{Extension information}

\begin{center}\begin{tabular}{c|cccccccccccccc}
Block 1 &$1$&$8$&$16$&$26$&$48$&$128$&$160$&$246$&$416$&$768$&$784$&$2560$&$3936$&$4096$
\\ \hline 
$1$&$0$&$1$&$0$&$1$&$1$&$0$&$0$&$1$&$1$&$0$&$0$&$0$&$0$&$0$
\\$8$&$1$&$0$&$0$&$1$&$0$&$0$&$0$&$0$&$0$&$0$&$1$&&&
\\$16$&$0$&$0$&$0$&$1$&$0$&$0$&$0$&$0$&$2$&$1$&&&&
\\$26$&$1$&$1$&$1$&$0$&$1$&$0$&$1$&$2$&&&&&&
\\$48$&$1$&$0$&$0$&$1$&$0$&$0$&$0$&&&&&&&
\\$128$&$0$&$0$&$0$&$0$&$0$&&&&&&&&&
\\$160$&$0$&$0$&$0$&$1$&$0$&&&&&&&&&
\\$246$&$1$&$0$&$0$&$2$&&&&&&&&&&
\\$416$&$1$&$0$&$2$&&&&&&&&&&&
\\$768$&$0$&$0$&$1$&&&&&&&&&&&
\\$784$&$0$&$1$&&&&&&&&&&&&
\end{tabular}\end{center}

\newpage
\section{\texorpdfstring{$\Sp_8(3)$}{Sp(8,3)}}

\subsection*{Block information}

These modules are only necessarily complete up to dimension 4000.

\begin{center}\begin{tabular}{ccc}
\hline Block & Kernel & Modules
\\\hline $1$ & $2$ & $1,27,36,41,266,279,513,594_1,594_2,784,1114,1891,4320,6237,6426,10150$, and 28 more
\\ $2$ & $2$ & $43046721$
\\ $3$ & $1$ & $8,40,112,288,504,1016,1072,1512,2072,2640,4752,6048,9080,$ and 23 more
\\\hline\end{tabular}\end{center}

\subsection*{Highest weight labels}

These are only up to dimension 5000.


\begin{center}\begin{tabular}{cccccccccccccc}
\hline $\lambda$&0000&0100&2000&0001&0200&1010&0101&2100&0002&0020&2001&1110&2200
\\\hline $L(\lambda)$ & $1$&$27$&$36$&$41$&$266$&$279$&$513$&$594_1$&$594_2$&$784$&$1114$&$1891$&$4320$
\\ $\dim(W(\lambda))$&$1$&$27$&$36$&$42$&$308$&$315$&$792$&$594$&$594$&$825$&$1155$&$4096$&$4914$
\\\hline
\end{tabular}\end{center}

\begin{center}\begin{tabular}{cccccccccccc}
\hline $\lambda$&1000&0010&1100&1001&0110&0011&2010&1200&1101&1002&1020
\\\hline $L(\lambda)$ & $8$&$40$&$112$&$288$&$504$&$1016$&$1072$&$1512$&$2072$&$2640$&$4752$
\\ $\dim(W(\lambda))$&$8$&$48$&$160$&$288$&$792$&$1056$&$1232$&$1512$&$4200$&$3696$&$4752$
\\\hline \end{tabular}\end{center}

\subsection*{Extension information}

\begin{center}\begin{tabular}{ccc}\hline
Module & Tested to & $\Ext^1=k$ 
\\\hline $1$ & $1891$ & $41$, $266$
\\ $27$ & $41$ & None
\\ $36$ & $41$ & None
\\ $41$ & $41$ & $1,41$
\\ \hline\end{tabular}\end{center}

\begin{center}\begin{tabular}{ccc}\hline
Module & Tested to & $\Ext^1=k$ 
\\\hline $8$ & $112$ & $40,112$
\\ $40$ & $112$ & $8,40$
\\ \hline\end{tabular}\end{center}

\newpage
\section{\texorpdfstring{$\Sp_8(5)$}{Sp(8,5)}}

\subsection*{Block information}

These modules are only necessarily complete up to dimension 4000.

\begin{center}\begin{tabular}{ccc}
\hline Block & Kernel & Modules
\\\hline $1$ & $2$ & $1,27,36,42,281,313,315,330,558,792,825,1155,1652,2183_1,2183_2,3138,3262$...
\\ $2$ & $2$ & $152587890625$
\\ $3$ & $1$ & $8,48,120,160,288,312,744,784,944,1504,2192,3144,3520$...
\\\hline\end{tabular}\end{center}

\subsection*{Highest weight labels}

These are only up to dimension 4000.

\begin{center}\begin{tabular}{cccccccccc}
\hline $\lambda$&0000&0100&2000&0001&0200&0002&1010&4000&2100
\\\hline $L(\lambda)$ & $1$&$27$&$36$&$42$&$281$&$313$&$315$&$330$&$558$
\\ $\dim(W(\lambda))$&$1$&$27$&$36$&$42$&$308$&$594$&$315$&$330$&$594$
\\\hline $\lambda$&0101&0020&2001&2200&0300&1011&3010&1110&
\\\hline $L(\lambda)$ &$792$&$825$&$1155$& $1652$&$2183_1$&$2183_2$&$3138$&$3262$&
\\ $\dim(W(\lambda))$&$792$&$825$&$1155$&$4914$&$2184$&$6237$&$3696$&$4096$&
\\\hline
\end{tabular}\end{center}

\begin{center}\begin{tabular}{cccccccccccccc}
\hline $\lambda$&1000&0010&3000&1100&1001&0011&0110&3100&2010&1200&1002&1300&3001
\\\hline $L(\lambda)$ & $8$&$48$&$120$&$160$&$288$&$312$&$744$&$784$&$944$&$1504$&$2192$&$3144$&$3520$
\\ $\dim(W(\lambda))$&$8$&$48$&$120$&$160$&$288$&$1056$&$792$&$1728$&$1232$&$1512$&$3696$&$9408$&$3520$
\\\hline
\end{tabular}\end{center}

\subsection*{Extension information}

\begin{center}\begin{tabular}{ccc}\hline
Module & Tested to & $\Ext^1=k$ 
\\\hline $1$ & $1155$ & None
\\ $27$ & $281$ & $281$
\\ $36$ & $42$ & None
\\ $42$ & $42$ & None
\\ \hline\end{tabular}\end{center}

\begin{center}\begin{tabular}{ccc}\hline
Module & Tested to & $\Ext^1=k$ 
\\\hline $8$ & $288$ & None
\\ $48$ & $48$ & None
\\ \hline\end{tabular}\end{center}

\newpage
\section{\texorpdfstring{$\Omega_7(3)$}{O7(3)}}

\subsection*{Block information}

\begin{center}\begin{tabular}{ccc}
\hline Block & Kernel & Modules
\\\hline Block 1& $2$ & $1$, $7$, $21$, $27$, $35$, $63$, $132$, $189_1$, $189_2$, $309$, $518$, $672_1$, $797$, $1602$, $2052$, $2855$, $5670$
\\ Block 2 & $2$ & $19683$
\\ Block 3 & $1$ & $8$, $48$, $104$, $120$, $384$, $672_2$, $1512$, $2400$, $7272$
\\\hline\end{tabular}\end{center}

\subsection*{Highest weight labels}

\begin{center}\begin{tabular}{cccccccccc}
\hline $\lambda$&000&100&010&200&002&110&020&102&012
\\\hline $L(\lambda)$ &$1$&$7$&$21$&$27$&$35$&$63$&$132$&$189_1$&$189_2$
\\ $\dim(W(\lambda))$&$1$&$7$&$21$&$27$&$35$&$105$&$168$&$189$&$378$
\\\hline $\lambda$&210&202&120&112&220&022&212&122&222
\\\hline $L(\lambda)$&$309$&$518$&$672_1$&$797$&$1602$&$2052$&$2855$&$5670$&$19683$
\\ $\dim(W(\lambda))$&$330$&$616$&$693$&$1617$&$1911$&$2079$&$4550$&$7722$&$19683$
\\\hline
\end{tabular}\end{center}

\begin{center}\begin{tabular}{cccccccccc}
\hline $\lambda$&001&101&011&201&111&021&211&121&221
\\\hline $L(\lambda)$ &$8$&$48$&$104$&$120$&$384$&$672_2$&$1512$&$2400$&$7272$
\\ $\dim(W(\lambda))$&$8$&$48$&$112$&$168$&$512$&$720$&$1512$&$2800$&$7392$
\\\hline
\end{tabular}\end{center}

\subsection*{Extension information}

\begin{small}\begin{center}\begin{tabular}{c|ccccccccccccccccc}
Block 1 & $1$&$7$&$21$&$27$&$35$&$63$&$132$&$189_1$&$189_2$&$309$&$518$&$672_1$&$797$&$1602$
&$2052$&$2855$&$5670$
\\ \hline $1$&$0$&$0$&$0$&$0$&$0$&$1$&$3$&$0$&$0$&$0$&$0$&$0$&$2$&$0$&$0$&$2$&$0$
\\ $7$&$0$&$0$&$0$&$0$&$0$&$3$&$1$&$0$&$1$&$1$&$0$&&&&&&$0$
\\ $21$&$0$&$0$&$0$&$0$&$0$&$1$&$2$&$0$&$0$&&&&&&&&$0$
\\ $27$&$0$&$0$&$0$&$0$&$0$&$0$&$1$&$0$&$0$&&&&&&&&$0$
\\ $35$&$0$&$0$&$0$&$0$&$0$&$2$&$2$&$0$&$0$&&&&&&&&$0$
\\ $63$&$1$&$3$&$1$&$0$&$2$&$0$&$0$&$0$&$0$&&&&&&&&$0$
\\ $132$&$3$&$1$&$2$&$1$&$2$&$0$&$0$&&&&&&&&&&$1$
\\ $189_1$&$0$&$0$&$0$&$0$&$0$&$0$&&&&&&&&&&&$0$
\\ $189_2$&$0$&$1$&$0$&$0$&$0$&$0$&&&&&&&&&&&$0$
\\ $309$&$0$&$1$&&&&&&&&&&&&&&&$0$
\\ $518$&$0$&$0$&&&&&&&&&&&&&&&$0$
\\ $672_1$&$0$&&&&&&&&&&&&&&&&$0$
\\ $797$&$2$&&&&&&&&&&&&&&&&$0$
\\ $1602$&$0$&&&&&&&&&&&&&&&&$0$
\\ $2052$&$0$&&&&&&&&&&&&&&&&$1$
\\ $2855$&$2$&&&&&&&&&&&&&&&&$0$
\\ $5670$&$0$&$0$&$0$&$0$&$0$&$0$&$1$&$0$&$0$&$0$&$0$&$0$&$0$&$0$&$1$&$0$&$0$
\end{tabular}\end{center}\end{small}


\begin{center}\begin{tabular}{c|cccccc}
Block 3&$8$&$48$&$104$&$120$&$384$&$672_2$
\\ \hline $8$&$0$&$0$&$2$&$1$&$6$&$2$
\\ $48$&$0$&$0$&$0$&$2$&&
\\ $104$&$2$&$0$&$0$&$0$&&
\\ $120$&$1$&$2$&$0$&&&
\\ $384$&$6$&&&&&
\\ $672_2$&$2$&&&&&
\end{tabular}\end{center}

\newpage
\section{\texorpdfstring{$\Omega_7(5)$}{O7(5)}}

\subsection*{Block information}

These modules are only necessarily complete up to dimension 2000.

\begin{center}\begin{tabular}{ccc}
\hline Block & Kernel & Modules
\\\hline Block 1& $2$ & $1$, $7$, $21$, $27$, $35$, $77$, $105$, $141$, $168_1$, $182$, $293$, $330$, $371_1$, $371_2$, $448_1$, $483$, $616$, $825$, $988$,
\\ && $1093$, $1246$, $1322$, $1513$, $1560$, $1750$, $2008$, $2310$, $2380$, $2393$, $3375$, $4095$,
$4550$, $4740$,
\\ && $4920$, $5622$, $6608$, $7372$, $7595$, $7783$, $7916$, $9625$, $10800$, $12125$, $12704$, $16275$, $19879$,
\\ && $20440$, $23625$, $26446$, $26684$, $29100$, $30472$, and 22 more.
\\ Block 2 & $2$ & $1953125$
\\ Block 3 & $1$ & $8$, $48$, $64$, $104$, $168_2$, $344$, $448_2$, $560$, $664$, $720$, $952$, $1008$, $1064$, $1624$, $2360$,
\\ &&  $2800$, $2912$, $2968$, $3752$,  $4200$, $4584$, $4600$, $4992$, $7112$, $8408$,
\\ && $8616$, $9416$, $14936$, $17080$, $28400$, and 20 more.
\\\hline\end{tabular}\end{center}

\subsection*{Highest weight labels}

These are only up to dimension 2000.

\begin{center}\begin{tabular}{cccccccccccccc}
\hline $\lambda$&000&100&010&200&002&300&110&020&102&400&004&210&012
\\\hline $L(\lambda)$ &$1$&$7$&$21$&$27$&$35$&$77$&$105$&$141$&$168_1$&$182$&$293$&$330$&$371_1$
\\ $\dim(W(\lambda))$&$1$&$7$&$21$&$27$&$35$&$77$&$105$&$168$&$189$&$182$&$294$&$330$&$378$
\\\hline $\lambda$&310&202&220&120&030&130&104&112&022&040&302&410&
\\\hline $L(\lambda)$&$371_2$&$448_1$& $483$&$616$&$825$&$988$&$1093$&$1246$&$1322$&$1513$&$1560$&$1750$&
\\ $\dim(W(\lambda))$&$819$&$616$&$1911$&$693$&$825$&$3003$&$1386$&$1617$&$2079$&$3003$&$1560$&$1750$&
\\\hline
\end{tabular}%
\end{center}

\subsection*{Extension information}


\begin{center}%
\end{center}

\newpage
\section{\texorpdfstring{$\Omega_7(7)$}{O7(7)}}

\subsection*{Block information}

These modules are only necessarily complete up to dimension 2000.

\begin{center}%
\end{center}

\subsection*{Highest weight labels}

\begin{center}%
\end{center}

\subsection*{Extension information}

%
%
%

\begin{center}

\end{center}

\newpage
\section{\texorpdfstring{$\Omega_9(3)$}{O9(3)}}

\subsection*{Block information}

These modules are only necessarily complete up to dimension 4000.

\begin{center}%
\end{center}

\subsection*{Highest weight labels}

These are only up to dimension 4000.

\begin{center}%
\end{center}

\subsection*{Extension information}

\begin{center}\begin{tabular}{ccc}\hline
Module & Tested to & $\Ext^1=k$
\\ \hline $1$ & $579$ & $43,369,579$
\\ $9$ & $84$ & None
\\ $36$ & $84$ & None
\\ $43$ & $84$ & $1$
\\ $84$ & $84$ & None
\\ \hline\end{tabular}\end{center}

\newpage
\section{\texorpdfstring{$\Omega_9(5)$}{O9(5)}}

\subsection*{Block information}

These modules are only necessarily complete up to dimension 4000.

\begin{center}%
\end{center}

\subsection*{Highest weight labels}

These are only up to dimension 4000.

\begin{center}%
\end{center}

\subsection*{Extension information}

\begin{center}%
\end{center}

\newpage
\section{\texorpdfstring{$\Omega_8^+(2)$}{O8+(2)}}

\subsection*{Block information}

\begin{center}%
\end{center}

\subsection*{Identification of modules}

Since there is an $S_3$ of automorphisms permuting the $8$-dimensional module, we can make any choice for $8_1$, $8_2$ and $8_3$. We then define $8_i\otimes 26=48_i\oplus 160_i$, and $748_i$ as the socle of $8_i\otimes 246$.

\subsection*{Highest weight labels}

\begin{center}%
\end{center}

\subsection*{Extension information}

\begin{center}%
\end{center}

\subsection*{Symmetric and exterior squares of modules}

\begin{center}%
\end{center}

\subsection*{Actions of unipotent elements on modules}

\begin{center}%
\end{center}

\subsection*{Traces of semisimple elements on modules}

\begin{small}\begin{center}%
\end{center}\end{small}

\newpage
\section{\texorpdfstring{$\Omega_8^+(3)$}{O8+(3)}}

\subsection*{Block information}

\begin{center}%
\end{center}

\subsection*{Identification of modules}

Since there is an $S_3$ of automorphisms permuting the $8$-dimensional module, we can make any choice for $8_1$, $8_2$ and $8_3$. Then $35_i$ lies in $S^2(8_i)$, $518_{i+2}$ lies in $35_i\otimes 35_{i+1}$ (with the indices taken modulo $3$), and $567_i$ lies in $28\otimes 35_i$. Finally, $3948_i$ is a composition factor of $35_i\otimes 195$.

\subsection*{Highest weight labels}

\noindent Blocks 1 and 2:
\begin{center}%
\end{center}

\subsection*{Extension information}

\begin{center}%
\end{center}

\newpage
\section{\texorpdfstring{$\Omega_8^+(5)$}{O8+(5)}}

\subsection*{Block information}
These modules are only necessarily complete up to dimension 10000.
\begin{center}%
\end{center}

\subsection*{Highest weight labels}

\begin{center}%
\end{center}

\subsection*{Extension information}

\begin{center}%
\end{center}

\subsection*{Actions of unipotent elements on modules}

\begin{center}%
\end{center}

\newpage
\section{\texorpdfstring{$\Omega_8^-(2)$}{O8-(2)}}

\subsection*{Block information}

\begin{center}%
\end{center}

\subsection*{Identification of modules}

Define $8_1$, $48_1$, $160_1$, $784_1$ as the modules defined over $\F_2$. For the modules $n_2$ and $n_3$ for the various $n$, we chose $8_2$ and $8_3$ arbitrarily, since they are permuted by the outer automorphism, and define the rest by $8_i\otimes 26=48_i\oplus 160_i$ and $784_i=\soc(8_i\otimes 246)$, as with $\Omega_8^+(2)$.

\subsection*{Extension information}

\begin{center}%
\end{center}


\subsection*{Actions of unipotent elements on modules}

\begin{center}%
\end{center}

\subsection*{Traces of semisimple elements on modules}

\begin{small}\begin{center}%
\end{center}\end{small}

\newpage
\section{\texorpdfstring{$\Omega_8^-(3)$}{O8-(3)}}

\subsection*{Block information}

\begin{center}%
\end{center}

\subsection*{Extension information}
\begin{center}%
\end{center}

\newpage
\section{\texorpdfstring{$\Omega_8^-(5)$}{O8-(5)}}

\subsection*{Block information}
\begin{center}%
\end{center}

\subsection*{Extension information}
\begin{center}%
\end{center}

\newpage
\section{\texorpdfstring{${}^3\!D_4(2)$}{3D4(2)}}

\subsection*{Block information}

\begin{center}%
\end{center}

\subsection*{Identification of modules}

Choose $8_1$ arbitrarily from the three $8_i$. Define $8_{i+1}$ to lie in the socle of $S^2(8_i)$. As with $\Omega_8^+(2)$, we then define $8_i\otimes 26=48_i\oplus 160_i$, and $748_i$ as the socle of $8_i\otimes 246$.

\subsection*{Extension information}

\begin{center}%
\end{center}

\newpage
\section{\texorpdfstring{${}^3\!D_4(3)$}{3D4(3)}}

\subsection*{Block information}

\begin{center}%
\end{center}

\subsection*{Extension information}

\begin{center}
%
\end{center}

\newpage
\section{\texorpdfstring{${}^3\!D_4(5)$}{3D4(5)}}

\subsection*{Block information}
These modules are only necessarily complete up to dimension 10000.
\begin{center}%
\end{center}

\subsection*{Extension information}
\begin{center}%
\end{center}

To distinguish $8_i$ and $104_i$, note that $104_i$ is a composition factor of $S^3(8_i)$, and in fact $S^3(8_i)=8_i/104_i/8_i$.

\newpage
\section{\texorpdfstring{$F_4(2)$}{F4(2)}}

\subsection*{Block information}

\begin{center}%
\end{center}

\subsection*{Identification of modules}

Define $26_1$ arbitrarily from the $26$-dimensionals. We define $246_i$ to lie in $\Lambda^2(26_i)$. The module $676=26_1\otimes 26_2$. We have that $4096_i$ lies in $26_i\otimes 246_i$. Then $6376_i=26_i\otimes 246_{3-i}$.

\subsection*{Highest weight labels}

\begin{center}%
\end{center}

\subsection*{Extension information}

\begin{center}%
\end{center}

\newpage
\section{\texorpdfstring{$G_2(2)'$}{G2(2)'}}

\subsection*{Block information}

\begin{center}%
\end{center}

\subsection*{Highest weight labels}

\begin{center}%
\end{center}

\subsection*{Extension information}

\begin{center}%
\end{center}

\subsection*{Socle layers of projective modules}

\[%
\]

\subsection*{Symmetric and exterior squares of modules}

\begin{center}%
\end{center}

\subsection*{Actions of unipotent elements on modules}

\begin{center}%
\end{center}

\newpage
\section{\texorpdfstring{$G_2(3)$}{G2(3)}}

\subsection*{Block information}

\begin{center}%
\end{center}

\subsection*{Highest weight labels}

\begin{center}%
\end{center}

\subsection*{Extension information}

\begin{center}%
\end{center}

\subsection*{Symmetric and exterior squares of modules}

\begin{center}%
\end{center}

\subsection*{Actions of unipotent elements on modules}

\begin{center}%
\end{center}

\newpage
\section{\texorpdfstring{$G_2(4)$}{G2(4)}}

\subsection*{Block information}

\begin{center}%
\end{center}

\subsection*{Identification of modules}

\begin{enumerate}
\item $84_{1,2}=6_1\otimes 14_2$, $84_{2,1}=14_1\otimes 6_2$
\item $384_{1,2}=6_1\otimes 64_2$, $384_{2,1}=64_1\otimes 6_2$
\item $896_{1,2}=14_1\otimes 64_2$, $896_{2,1}=64_1\otimes 14_2$
\end{enumerate}

\subsection*{Extension information}

\begin{center}%
\end{center}

\subsection*{Symmetric and exterior squares of modules}

\begin{center}%
\end{center}

\subsection*{Actions of unipotent elements on modules}

\begin{center}%
\end{center}

\newpage
\section{\texorpdfstring{$G_2(5)$}{G2(5)}}

\subsection*{Block information}

\begin{center}%
\end{center}

\subsection*{Highest weight labels}

\begin{center}%
\end{center}

\subsection*{Extension information}

There are no extensions between $8456$ and modules of dimension greater than $792$.

\begin{center}%
\end{center}

\newpage
\section{\texorpdfstring{$G_2(7)$}{G2(7)}}

\subsection*{Block information}

\begin{center}%
\end{center}

\subsection*{Highest weight labels}

\begin{center}%
\end{center}

\subsection*{Extension information}

\begin{center}%
\end{center}

\subsection*{Actions of unipotent elements on modules}

\begin{center}%
\end{center}

\subsection*{Traces of semisimple elements on modules}

\begin{center}%
\end{center}

\end{document}